\documentclass{amsart}%
\usepackage{graphicx}
\usepackage{amssymb}%
\newtheorem{theorem}{Theorem}[section]
\newtheorem{corollary}[theorem]{Corollary}

\newtheorem{lemma}[theorem]{Lemma}
\newtheorem{proposition}[theorem]{Proposition}

\theoremstyle{definition}
\newtheorem{definition}[theorem]{Definition}

\newtheorem{examples}[theorem]{Examples}

\newtheorem{remark}[theorem]{Remark}
\newtheorem{final remarks and questions}[theorem]{Final remarks and questions}

\numberwithin{equation}{section}
\begin{document}

\title{On the Equations Defining Toric L.C.I.-Singularities}
\author{Dimitrios I. Dais}
\address{Department of Mathematics and Statistics,
         University of Cyprus, P.O. Box 20537, CY-1678 Nicosia, Cyprus}
\email{ddais@ucy.ac.cy}
\author{Martin Henk}
\address{Technische Universit\"{a}t Wien, Abteilung f\"{u}r
Analysis, Wiedner Hauptstr. 8-10, A-1040 Wien, Austria}
\curraddr{Technische Universit{\"a}t Berlin, Institut f{\"u}r Mathematik, Sekr. MA 6-2, Stra{\ss}e des 17. Juni 136, D-10623 Berlin, Germany}
\email{henk@math.tu-berlin.de}
\thanks{The second author would like to
thank the Mathematics Department of the University of Crete for hospitality
and support during the spring term 2001, where this work was initiated.}
\date{}

\begin{abstract}
\noindent Based on Nakajima's Classification Theorem \cite{Nakajima} we
describe the precise form of the binomial equations which determine toric
locally complete intersection (``l.c.i'') singularities.
\end{abstract}

\maketitle

\section{ Introduction}

\noindent{}An affine toric variety $U_{\sigma}=$ Spec$\left(
\mathbb{C}\left[ M\cap\sigma^{\vee}\right] \right) $ associated to
a rational strongly convex polyhedral cone $\sigma$ (where
rank$(M)=$ dim$(\sigma)=d$) is Gorenstein if and only if $\sigma$
supports a $(d-1)$-dimensional lattice polytope $P$ (w.r.t. $N=$
Hom$_{\mathbb{Z}}\left( M,\mathbb{Z}\right) $) lying on a
``primitive'' affine hyperplane of the form $\left\{ \mathbf{x}\in N_{%
\mathbb{R}}\ \left| \ \left\langle
m_{\sigma},\mathbf{x}\right\rangle =1\right. \right\} $ (up to
lattice automorphism). Nakajima \cite{Nakajima} classified in 1986
all affine toric locally complete intersection (``l.c.i.'')
varieties $U_{\sigma}$ by providing a suitable parametrization for
all the corresponding polytopes $P.$ More recently, the class of
toric l.c.i.-singularities turned out to have some properties of
particular importance in both algebraic and geometric aspects. For
instance,\smallskip\

\noindent{}(i) the algebras $\mathbb{C}\left[ N\cap\sigma\right] $
have the Koszul property (cf. \cite{BGT}, \cite{DHaZ}),\smallskip

\noindent{}(ii) all set-theoretic complete intersections of
binomial hypersurfaces in an affine complex space are affine toric
ideal-theoretic complete intersections (also in a more general
sense, including the non-normal ones, cf. \cite[Thm.
4]{BMT}),\smallskip

\noindent{}(iii) all toric l.c.i.-singularities admit projective
crepant resolutions (see \cite{DHZ}, \cite{DHaZ}), and\smallskip

\noindent{}(iv) the $i$-th jet schemes of the underlying spaces
$U_{\sigma}$
of toric l.c.i.-singularities are irreducible for all $i\geq1$ (see \cite[%
Thm. 4.13]{Mustata}).\smallskip

The Main Theorem of the present paper (see below Thm. \ref{MAIN})
gives a
precise description of the binomial-type equations defining singular l.c.i. $%
U_{\sigma}$'s in terms of the corresponding admissible
free-parameter
sequence (or ``matrix'') $\mathbf{m}$ of the Nakajima polytope $P\sim P_{%
\mathbf{m}}^{\left( d\right) }$, and generalizes a result of Ishida \cite[\S %
8]{Ishida} which concerns dimensions $2$ and $3.\medskip$

\noindent In this section, we introduce the brief ``toric glossary'' \textsf{%
(a)}-\textsf{(e) }and the notation which will be used in the
sequel. For further details on the theory of toric geometry the
reader is referred to the textbooks of Oda \cite{Oda}, Fulton
\cite{Fulton} and Ewald \cite{Ewald}, and to the lecture notes
\cite{KKMS}. \medskip

\noindent\textsf{(a)} The \textit{linear hull, }the\textit{\
affine hull}, the \textit{positive hull} and \textit{the convex
hull} of a set $B$ of vectors of $\mathbb{R}^{r}$, $r\geq1,$ will
be denoted by lin$\left(
B\right) $, aff$\left( B\right) $, pos$\left( B\right) $ (or $\mathbb{R}%
_{\geq0}\,B$) and conv$\left( B\right) $, respectively. The \textit{dimension%
} dim$\left( B\right) $ of a $B\subset\mathbb{R}^{r}$ is defined
to be the dimension of its affine hull. \newline
\newline
\textsf{(b) }Let $N$ be a free $\mathbb{Z}$-module of rank
$r\geq1$. $\ N$
can be regarded as a \textit{lattice }within \ $N_{\mathbb{R}}:=N\otimes_{%
\mathbb{Z}}\mathbb{R}\cong\mathbb{R}^{r}$. The \textit{lattice
determinant} det$\left( N\right) $ of $N$ is the $r$-volume of the
parallelepiped spanned by any $\mathbb{Z}$-basis of it. An $n\in
N$ is called \textit{primitive} if conv$\left( \left\{
\mathbf{0},n\right\} \right) \cap N$ contains no other points
except $\mathbf{0}$ and $n$.\smallskip

Let $N$ be as above, $M:=$ Hom$_{\mathbb{Z}}\left(
N,\mathbb{Z}\right) $ its dual lattice,
$N_{\mathbb{R}},M_{\mathbb{R}}$ their real scalar extensions,
and $\left\langle .,.\right\rangle :M_{\mathbb{R}}\times N_{\mathbb{R}%
}\rightarrow \mathbb{R}$ the natural $\mathbb{R}$-bilinear
pairing. A subset
$\sigma $ of $N_{\mathbb{R}}$ is called \textit{convex polyhedral cone} (%
\textit{c.p.c.}, for short) if there exist $n_{1},\ldots ,n_{k}\in N_{%
\mathbb{R}}$, such that $\sigma =$ pos$\left( \left\{ n_{1},\ldots
,n_{k}\right\} \right) $. Its \textit{relative interior
}int$\left( \sigma \right) $ is the usual topological interior of
it, considered as subset of
lin$\left( \sigma \right) =\sigma +\left( -\sigma \right) $. The \textit{%
dual cone} $\sigma ^{\vee }$ of a c.p.c. $\sigma $ is a c.p. cone
defined by
\begin{equation*}
\sigma ^{\vee }:=\left\{ \mathbf{y}\in M_{\mathbb{R}}\ \left\vert
\
\left\langle \mathbf{y},\mathbf{x}\right\rangle \geq 0,\ \forall \mathbf{x}%
,\ \mathbf{x}\in \sigma \right. \right\} \;.\;
\end{equation*}%
Note that $\left( \sigma ^{\vee }\right) ^{\vee }=\sigma $ and
\begin{equation*}
\text{dim}\left( \sigma \cap \left( -\sigma \right) \right) +\text{dim}%
\left( \sigma ^{\vee }\right) =\text{dim}\left( \sigma ^{\vee
}\cap \left( -\sigma ^{\vee }\right) \right) +\text{dim}\left(
\sigma \right) =r.
\end{equation*}%
A subset $\tau $ of a c.p.c. $\sigma $ is called a \textit{face}
of $\sigma $ (notation: $\tau \prec \sigma $), if for some
$m_{0}\in \sigma ^{\vee }$ we
have $\tau =\left\{ \mathbf{x}\in \sigma \ \left\vert \ \left\langle m_{0},%
\mathbf{x}\right\rangle =0\right. \right\} $. In particular,
1-dimensional faces are called \textit{rays}.

A c.p.c. $\sigma=$ pos$\left( \left\{ n_{1},\ldots,n_{k}\right\}
\right) $ is called \textit{simplicial} (resp. \textit{rational})
if $n_{1},\ldots
,n_{k}$ are $\mathbb{R}$-linearly independent (resp. if $n_{1},\ldots,n_{k}%
\in N_{\mathbb{Q}}$, where $N_{\mathbb{Q}}:=N\otimes_{\mathbb{Z}}\mathbb{Q}$%
). If $\varrho$ is a ray of a rational c.p.c. $\sigma$, then we denote by $%
n\left( \varrho\right) \in N\cap\varrho$ \ the unique primitive vector with $%
\varrho=\mathbb{R}_{\geq0}\ n\left( \varrho\right) ,$ and we set
\begin{equation*}
\text{Gen}\left( \sigma\right) :=\left\{ n\left( \varrho\right) \
\left\vert \ \varrho\text{ rays of }\sigma\right. \right\} .
\end{equation*}
A \textit{strongly convex polyhedral cone }(\textit{s.c.p.c.}, for
short) is a c.p.c. $\sigma$ for which $\sigma\cap\left(
-\sigma\right) =\left\{ \mathbf{0}\right\} $, i.e., for which
dim$\left( \sigma^{\vee}\right) =r$.
The s.c.p. cones are alternatively called \textit{pointed cones} (having $%
\mathbf{0}$ as their apex).\newline
\newline
\textsf{(c) }If $\sigma\subset N_{\mathbb{R}}$ is a rational
s.c.p.c., then the subsemigroup $\sigma\cap N$ of $N$ is a monoid.
The following proposition follows from results due to Gordan,
Hilbert and van der Corput (cf. \cite[Thm. 16.4, p.
233]{Schrijver}).

\begin{proposition}[Minimal generating system]
\label{MINGS}$\sigma\cap N$ is finitely generated as additive
semigroup for every rational c.p.c. $\sigma\subset
N_{\mathbb{R}}$. Moreover, if $\sigma$ is strongly convex, then
among all the systems of generators of $\sigma\cap N
$, there is a system $\mathbf{Hlb}_{N}\left( \sigma\right) $ of \emph{%
minimal cardinality}, which is uniquely determined \emph{(}up to
the
ordering of its elements\emph{)} by the following characterization\emph{%
:\smallskip}
\begin{equation}
\mathbf{Hilb}_{N}\left( \sigma\right) =\left\{
n\in\sigma\cap\left( N\smallsetminus\left\{ \mathbf{0}\right\}
\right) \ \left\vert \
\begin{array}{l}
n\ \text{\emph{cannot be expressed }} \\
\text{\emph{as the sum of two other }} \\
\text{\emph{vectors belonging }} \\
\text{\emph{to\ } }\sigma\cap\left( N\smallsetminus\left\{ \mathbf{0}%
\right\} \right)%
\end{array}
\right. \right\} .   \label{Hilbbasis}
\end{equation}
$\mathbf{Hilb}_{N}\left( \sigma\right) $ \emph{is called
}\textit{the Hilbert basis of }$\sigma$ w.r.t. $N.$
\end{proposition}

\noindent

\noindent\textsf{(d)} For a lattice $N$ of rank $r$ having $M$ as
its dual, we define an $r$-dimensional \textit{algebraic torus
}$T_{N}\cong\left(
\mathbb{C}^{\ast}\right) ^{r}$ by setting $T_{N}:=$Hom$_{\mathbb{Z}}\left( M,%
\mathbb{C}^{\ast}\right) =N\otimes_{\mathbb{Z}}\mathbb{C}^{\ast}.$ Every $%
m\in M$ assigns a character $\mathbf{e}\left( m\right) :T_{N}\rightarrow%
\mathbb{C}^{\ast}$. Moreover, each $n\in N$ determines an
$1$-parameter subgroup
\begin{equation*}
\vartheta_{n}:\mathbb{C}^{\ast}\rightarrow T_{N}\ \ \ \text{with\ \ \ }%
\vartheta_{n}\left( \lambda\right) \left( m\right)
:=\lambda^{\left\langle m,n\right\rangle }\text{, \ \ for\ \ \
}\lambda\in\mathbb{C}^{\ast},\ m\in M\ .\
\end{equation*}
We can therefore identify $M$ with the character group of $T_{N}$
and $N$ with the group of $1$-parameter subgroups of $T_{N}$. On
the other hand, for a rational s.c.p. cone $\sigma$ with
\begin{equation*}
M\cap\sigma^{\vee}=\mathbb{Z}_{\geq0}\
m_{1}+\cdots+\mathbb{Z}_{\geq0}\ m_{\nu},
\end{equation*}
we associate to the finitely generated monoidal subalgebra
\begin{equation*}
\mathbb{C}\left[ M\cap\sigma^{\vee}\right] =\bigoplus\limits_{m\in
M\cap\sigma^{\vee}}\mathbf{e}\left( m\right)
\end{equation*}
of the $\mathbb{C}$-algebra $\mathbb{C}\left[ M\right] =\bigoplus
\limits_{m\in M}\mathbf{e}\left( m\right) $ a \textit{toric affine
variety}
\begin{equation*}
U_{\sigma}:=U_{\sigma,N}:=\text{Spec}\left( \mathbb{C}\left[ M\cap
\sigma^{\vee}\right] \right) ,
\end{equation*}
which can be identified with the set of semigroup homomorphisms :
\begin{equation*}
U_{\sigma}=\left\{ u:M\cap\sigma^{\vee}\ \rightarrow\mathbb{C\
}\left\vert
\begin{array}{c}
\ u\left( \mathbf{0}\right) =1,\ u\left( m+m^{\prime}\right)
=u\left(
m\right) \cdot u\left( m^{\prime}\right) ,\smallskip\  \\
\text{for all \ \ }m,m^{\prime}\in M\cap\sigma^{\vee}%
\end{array}
\right. \right\} \ ,
\end{equation*}
where $\mathbf{e}\left( m\right) \left( u\right) :=u\left(
m\right) ,\
\forall m,\ m\in M\cap\sigma^{\vee}\ $ and\ $\forall u,\ u\in U_{\sigma}$%
.\medskip

\noindent {}$U_{\sigma }$ admits a canonical $T_{N}$-action which
extends
the group multiplication of the algebraic torus $T_{N}=U_{\left\{ \mathbf{0}%
\right\} }$:
\begin{equation}
T_{N}\times U_{\sigma }\ni \left( t,u\right) \longmapsto t\cdot
u\in U_{\sigma }  \label{torus action}
\end{equation}%
where, for $u\in U_{\sigma }$, $\left( t\cdot u\right) \left(
m\right) :=t\left( m\right) \cdot u\left( m\right) ,\ \forall m,\
m\in M\cap \sigma ^{\vee }$ . The orbits w.r.t. the action
(\ref{torus action}) are parametrized by the set of all the faces
of $\sigma $. For a $\tau \prec \sigma $, we denote by orb$\left(
\tau \right) $ the orbit which is associated to $\tau $.

\begin{proposition}[Embedding by binomials]
\label{EMB}In the analytic category, $U_{\sigma}$, identified with
its image under the injective map
\begin{equation*}
\left( \mathbf{e}\left( m_{1}\right) ,\ldots,\mathbf{e}\left(
m_{\nu }\right) \right)
:U_{\sigma}\hookrightarrow\mathbb{C}^{\nu},
\end{equation*}
can be regarded as an analytic set determined by a finite number
of equations of the form\emph{:} \emph{(monomial) = (monomial).}
This analytic structure induced on $U_{\sigma}$ is independent of
the semigroup generators $\left\{ m_{1},\ldots,m_{\nu}\right\} $
and each map $\mathbf{e}\left(
m\right) $ on $U_{\sigma}$ is holomorphic w.r.t. it. In particular, for $%
\tau\prec\sigma$, $U_{\tau}$ is an open subset of $U_{\sigma}$.
Moreover, if \emph{dim}$\left( \sigma\right) =r$ and
\begin{equation*}
\#\left( \mathbf{Hilb}_{M}\left( \sigma^{\vee}\right) \right) =k\
\ \left( \leq\nu\right) ,
\end{equation*}
then \emph{(}by Prop. \emph{\ref{MINGS})} $k$ is nothing but the \emph{%
(minimal)} \emph{embedding dimension} of $U_{\sigma}$, i.e. the \emph{minimal%
} number of generators of the maximal ideal of the local $\mathbb{C}$%
-algebra $\mathcal{O}_{U_{\sigma},\mathbf{0}}$.
\end{proposition}

\noindent See \cite[Prop. 1.2 and 1.3., pp. 4-7]{Oda} for a proof, and \cite%
{BRS}, \cite{BR}, \cite{EIS-STU}, \cite{H-Sh}, \cite{Sturmfels1}, \cite%
{Sturmfels2} for the general theory of the defining equations of
toric varieties.

\begin{remark}
\label{LATIDEAL}In fact, Prop. \ref{EMB} informs us that
\begin{equation*}
\mathbb{C}\left[ M\cap\sigma^{\vee}\right] \cong\mathbb{C}\left[
z_{1},\ldots,z_{\nu}\right] \,/\,\mathcal{I},
\end{equation*}
where the prime ideal $\mathcal{I}$ is generated by binomials
belonging to
\begin{equation*}
\left\{
\prod\limits_{i=1}^{\nu}z_{i}^{\kappa_{i}}-\prod\limits_{i=1}^{\nu
}z_{i}^{\xi_{i}}\,\,\left|
\begin{array}{l}
\left( \kappa_{1},\ldots,\kappa_{\nu}\right) ,\left(
\xi_{1},\ldots ,\xi_{\nu}\right) \in\left(
\mathbb{Z}_{\geq0}\right) ^{\nu}\smallskip \text{
} \\
\text{and }\left(
\kappa_{1}-\xi_{1},\ldots,\kappa_{\nu}-\xi_{\nu}\right)
\in L%
\end{array}
\right. \right\} ,
\end{equation*}
and
\begin{equation*}
L=\left\{ \left( \ell_{1},\ldots,\ell_{\nu}\right)
\in\mathbb{Z}^{\nu }\ \left| \
\sum_{i=1}^{\nu}\,\ell_{i}\,m_{i}=0\right. \right\}
\end{equation*}
(under an appropriate identification $M\cong\mathbb{Z}^{r}$). Let us fix a $%
\mathbb{Z}$-basis $\{\upsilon_{1},\ldots,\upsilon_{k}\}$ of the
integer
lattice $L$ and denote by $\mathcal{B}$ the $\left( \nu\times k\right) $%
-matrix having $\upsilon_{1},\ldots,\upsilon_{k}$ as its
column-vectors. Then the ideal
\begin{equation*}
\mathcal{J}_{\mathcal{B}}:=\left(
f_{\upsilon_{1}},f_{\upsilon_{2}},\ldots,f_{\upsilon_{k}}\right) \subset%
\mathbb{C}\left[ z_{1},\ldots,z_{\nu }\right] ,
\end{equation*}
generated by the binomials
\begin{equation*}
f_{\upsilon_{i}}:=\prod\limits_{j=1}^{\nu}z_{i}^{\left(
\upsilon_{i}\right)
_{j}^{+}}-\prod\limits_{j=1}^{\nu}z_{i}^{\left(
\upsilon_{i}\right) _{j}^{-}},\ \ 1\leq i\leq k,
\end{equation*}
with $\upsilon_{i}=\left( \left( \upsilon_{i}\right)
_{1},\ldots,\left(
\upsilon_{i}\right) _{\nu}\right) ^{\intercal},$ and%
\begin{equation*}
\left( \upsilon_{i}\right) _{j}^{+}:=\max\{0,\left(
\upsilon_{i}\right) _{j}\},\ \ \left( \upsilon_{i}\right)
_{j}^{-}:=\max\{0,-\left( \upsilon_{i}\right) _{j}\},\ 1\leq i\leq
k,\ 1\leq j\leq\nu,
\end{equation*}
is called the \textit{lattice ideal} associated to $\mathcal{B}$.
In general, we have
$\mathcal{J}_{\mathcal{B}}\subseteq\mathcal{I}$, with the
inclusion possibly strict.
\end{remark}

\begin{definition}[Dominating matrices]
A $\left( \nu\times k\right) $-matrix (with integer entries) is
called a \textit{mixed matrix} if every column of it has both a
positive and a negative entry. A mixed $\left( \nu\times k\right)
$-matrix is said to be \textit{dominating} if it does not contain
any mixed $\left( \rho\times\rho\right) $-submatrices for
$1\leq\rho\leq\min\{\nu,k\}.$
\end{definition}

\begin{theorem}
\label{SATUR}\emph{(i)}
$\mathcal{J}_{\mathcal{B}}\,\mathbb{C}\left[ z_{1}^{\pm 1},\ldots
,z_{\nu }^{\pm 1}\right] =\mathcal{I}\,\mathbb{C}\left[ z_{1}^{\pm
1},\ldots ,z_{\nu }^{\pm 1}\right] ,$ i.e.,
\begin{equation*}
\mathcal{J}_{\mathcal{B}}:(\prod\limits_{i=1}^{\nu }z_{i})^{\infty }=%
\mathcal{I}.\smallskip \newline
\end{equation*}%
\emph{(ii)} $\mathcal{J}_{\mathcal{B}}=\mathcal{I}$ if and only if $\mathcal{%
B}$ is a dominating matrix.
\end{theorem}

\noindent{}\textit{Proof}. For (i) we refer to \cite[Cor.
2.5]{EIS-STU} and \cite[Thm. 2.10]{BR}, and for (ii) to
\cite[Lemma 2.2]{F-Sh} or \cite[Thm.
1.1]{H-Sh}.\hfill$\square\bigskip$

\noindent\textsf{(e)}\textit{\ }The well-known hierarchy of
Noetherian rings:
\begin{equation*}
\text{(regular)}\Longrightarrow\text{(l.c.i.)}\Longrightarrow \text{%
(Gorenstein)}\Longrightarrow\text{(Cohen-Macaulay)}
\end{equation*}
(cf. \cite{Kunz}) is used to describe the punctual algebraic
behaviour of affine toric varieties.\textit{\ }

\begin{theorem}[Normality and CM-property]
The toric varieties $U_{\sigma}$ are always normal and
Cohen-Macaulay.
\end{theorem}

\noindent\textit{Proof. }For a proof of the normality property see \cite[%
Thm. 1.4, p. 7]{Oda}. The CM-property for toric varieties was
first shown by
Hochster in \cite{Hochster}. See also \cite[Thm. 14, p. 52]{KKMS}, and \cite[%
3.9, p. 125]{Oda}.\hfill$\square$\newline

\begin{definition}[Multiplicities and basic cones]
Let $N$\emph{\ }be a free\emph{\ }$\mathbb{Z}$-module of rank\emph{\ }$r$%
\emph{\ }and $\sigma \subset N_{\mathbb{R}}$ a\emph{\ }simplicial,
rational s.c.p.c. of dimension $d\leq r$. The cone $\sigma $ can
be obviously written
as $\sigma =\varrho _{1}+\cdots +\varrho _{d}$\emph{, }for distinct\emph{\ }%
rays $\varrho _{1},\ldots ,\varrho _{d}$. The \textit{multiplicity} mult$%
\left( \sigma ;N\right) $ of $\sigma $ with respect to $N$ is
defined as
\begin{equation*}
\text{mult}\left( \sigma ;N\right) :=\frac{\text{det}\left( \mathbb{Z}%
\,n\left( \varrho _{1}\right) \oplus \cdots \oplus
\mathbb{Z}\,n\left( \varrho _{d}\right) \right) }{\text{det}\left(
N_{\sigma }\right) },
\end{equation*}%
where $N_{\sigma }$ is the sublattice of $N$ generated (as
subgroup) by the set $N\cap $ lin$\left( \mathbb{\sigma }\right)
.$ \ If mult$\left( \sigma
;N\right) =1$,\emph{\ }then\emph{\ }$\sigma $\emph{\ }is called a \textit{%
basic cone} w.r.t. $N$.
\end{definition}

\begin{theorem}[Smoothness criterion]
\label{SMCR}The affine toric variety $U_{\sigma }=U_{\sigma,N}$ is
smooth
\emph{(}i.e., the corresponding local rings $\mathcal{O}_{U_{\sigma},\text{ }%
u}$ are regular at all points $u$ of $U_{\sigma}$\emph{)} iff
$\sigma$ is basic \textit{w.r.t.} $N$.
\end{theorem}

\noindent\textit{Proof. }See \cite[Ch. I, Thm. 4, p. 14]{KKMS}, and \cite[%
Thm. 1.10, p. 15]{Oda}.\hfill$\square\bigskip$\newline
Next Theorem describes necessary and sufficient conditions under which $%
U_{\sigma}$ is Gorenstein (cf. \cite[\S 7]{Ishida}).

\begin{theorem}[Gorenstein property]
\label{GOR-PR}Let $N$\emph{\ }be a free\emph{\ }$\mathbb{Z}$-module of rank%
\emph{\ }$r$\emph{\ }and $\sigma$ a\emph{\ }rational s.c.p. cone in $N_{%
\mathbb{R}}$ of dimension $d\leq r$. Then the following conditions
are equivalent\emph{:\medskip}\newline \emph{(i)} $\ \ U_{\sigma}$
is Gorenstein.\medskip\ \newline
\emph{(ii)} \ There exists an element $m_{\sigma}$ of $M$, such that $%
M\cap\left( \text{\emph{int}}\left( \sigma^{\vee}\right) \right)
=m_{\sigma}+(M\cap\sigma^{\vee}).$ \medskip\newline Moreover, if
$d=r$, then $m_{\sigma}$ in \emph{(ii) }is the unique primitive
element of $M\cap\left( \text{\emph{int}}\left(
\sigma^{\vee}\right) \right) $ with this property and the above
conditions are equivalent to the following
one\emph{:}\medskip\newline
\emph{(iii) Gen}$\left( \sigma\right) \subset\mathbf{H}$, where $\mathbf{H}%
=\left\{ \mathbf{x}\in N_{\mathbb{R}}\ \left| \ \left\langle m_{\sigma},%
\mathbf{x}\right\rangle =1\right. \right\} .$
\end{theorem}\bigskip

\section{Nakajima's Polytopes and Classification Theorem}

\noindent{}We shall henceforth focus our attention to Gorenstein
toric singularities and, in particular, to those which are locally
complete intersections (l.c.i.'s).\bigskip

\noindent \textsf{(a) }Let $N$ a free $\mathbb{Z}$-module of rank
$r\geq 2$ and $\sigma \subset N_{\mathbb{R}}$ a rational s.c.p.c.
of dimension $d\leq r $. Since $N/N_{\sigma }$ is torsion free,
there exists a lattice decomposition $N=N_{\sigma }\oplus
\breve{N}$, inducing a decomposition of
its dual $M=M_{\sigma }\oplus \breve{M}$, where $M_{\sigma }=$ Hom$_{\mathbb{%
Z}}\left( N_{\sigma },\mathbb{Z}\right) $ and $\breve{M}=$ Hom$_{\mathbb{Z}}(%
\breve{N},\mathbb{Z)}$. Writing $\sigma $ as $\sigma =\sigma
^{\prime
}\oplus \left\{ \mathbf{0}\right\} $ with $\sigma ^{\prime }$ a $d$%
-dimensional cone in $\left( N_{\sigma }\right) _{\mathbb{R}}$, we
obtain decompositions
\begin{equation*}
T_{N}\ \cong T_{N_{\sigma }}\ \times T_{\breve{N}}\ \ \ \ \
\text{and\ \ \ \ }M\cap \sigma ^{\vee }=\left( M_{\sigma }\cap
\left( \sigma ^{\prime }\right) ^{\vee }\right) \oplus \breve{M}\
,
\end{equation*}%
which give rise to the analytic isomorphisms:
\begin{equation*}
\begin{array}[b]{ccc}
U_{\sigma } & \cong \ U_{\sigma ^{\prime }}\times T_{\breve{N}}\
\cong \ U_{\sigma ^{\prime }}\times T_{N/N_{\sigma }}\ \cong  &
U_{\sigma ^{\prime
}}\times \left( \mathbb{C}^{\ast }\right) ^{r-d}\smallskip .%
\end{array}%
\end{equation*}%
$U_{\sigma }=U_{\sigma ,N}$ can be therefore viewed as a product of $%
U_{\sigma ^{\prime }}=U_{\sigma ^{\prime },N}$ \ and an $\left( r-d\right) $%
-dimensional algebraic torus. Obviously, the study of the
algebraic properties for $U_{\sigma }$ can be reduced to that of
the corresponding properties of $\ U_{\sigma ^{\prime }}$. For
instance, the \textit{singular locus} Sing$\left( U_{\sigma
}\right) $ of $U_{\sigma }$ equals
\begin{equation*}
\text{Sing}\left( U_{\sigma }\right) =\text{Sing}\left( U_{\sigma
^{\prime }}\right) \times \left( \mathbb{C}^{\ast }\right) ^{r-d}.
\end{equation*}%
In fact, the main reason for preferring to work with the affine variety $%
U_{\sigma ^{\prime }}$ (or with the germ $\left( U_{\sigma ^{\prime }},\text{%
orb}\left( \sigma ^{\prime }\right) \right) )$ instead of
$U_{\sigma }$, is
that since lin$\left( \sigma ^{\prime }\right) =\left( N_{\sigma }\right) _{%
\mathbb{R}}$, the orbit orb$\left( \sigma ^{\prime }\right) \in
U_{\sigma
^{\prime }}$ is the unique fixed closed point under the action of $%
T_{N_{\sigma }}$ on $U_{\sigma ^{\prime }}$.

\begin{definition}
\label{SINR}If $\sigma$\ is non-basic w.r.t. $N$, then
$U_{\sigma^{\prime}}$
will be called\emph{\ }\textit{the singular representative}\emph{\ }of $%
U_{\sigma}$ and orb$\left( \sigma^{\prime}\right) \in
U_{\sigma^{\prime}}$ the associated\emph{\
}\textit{distinguished}\emph{\ }singular point within
the\emph{\ }singular locus Sing$\left( U_{\sigma^{\prime}}\right) $ of $%
U_{\sigma^{\prime}}$.\emph{\ }
\end{definition}

\begin{definition}
\label{SPL}If $\sigma$\ is non-basic w.r.t. $N$, then it is also
useful to introduce the notion of the\emph{\
}\textquotedblleft\textit{splitting
codimension}\textquotedblright\ of orb$\left(
\sigma^{\prime}\right) \in U_{\sigma^{\prime}}$ as the number
\begin{equation*}
\text{max\emph{\ }}\left\{ \varkappa\in\left\{ 2,\ldots,d\right\}
\ \left\vert \
\begin{array}{l}
U_{\sigma^{\prime}}\cong U_{\sigma^{\prime\prime}}\times\mathbb{C}%
^{d-\varkappa}, \\
\text{for\ some\emph{\ }
}\sigma^{\prime\prime}\prec\sigma^{\prime}\text{ }
\\
\text{with \ dim}\left( \sigma^{\prime\prime}\right)
=\varkappa\text{ \ \ }
\\
\text{and \ Sing}\left( U_{\sigma^{\prime\prime}}\right) \neq\varnothing%
\end{array}
\right. \right\} \ .
\end{equation*}
If this number equals $d$, then $\left(
U_{\sigma^{\prime}},\text{orb}\left(
\sigma^{\prime}\right) \right) $\emph{\ }will be called an \textit{%
msc-singularity}, i.e., a singularity having the\emph{\ }maximum
splitting codimension.\smallskip
\end{definition}

\noindent\textsf{(b) }Gorenstein toric affine varieties are
completely determined by suitable lattice polytopes.

\begin{definition}[Lattice equivalence]
\label{LP}\emph{\ }If $N_{1}$ and $N_{2}$\ are two free
$\mathbb{Z}$-modules (not necessarily of\emph{\ }the same rank)
and $P_{1}\subset \left(
N_{1}\right) _{\mathbb{R}}$,\emph{\ }$P_{2}\subset \left( N_{2}\right) _{%
\mathbb{R}}$\emph{\ }two lattice polytopes with respect to them,
we shall say that $P_{1}$\ and $P_{2}$\emph{\ }are \textit{lattice
equivalent }to each other, and denote this by\emph{\ }$P_{1}\sim
P_{2}$, if $P_{1}$\ is
affinely equivalent to $P_{2}$\ via an\emph{\ }affine map\emph{\  }%
\begin{equation*}
\varpi :\left( N_{1}\right) _{\mathbb{R}}\rightarrow \left( N_{2}\right) _{%
\mathbb{R}}\emph{,}
\end{equation*}%
such that the restriction
\begin{equation*}
\varpi \left\vert _{\text{aff}\left( P\right) }\right.
:\text{aff}\left( P\right) \rightarrow \text{aff}\left( P^{\prime
}\newline \right)
\end{equation*}%
is a bijection mapping the polytope \emph{\ }$P_{1}$\emph{\ }onto
the (necessarily
equidimensional) polytope \emph{\ }$P_{2}$, every $i$-dimensional face of%
\emph{\ }$P_{1}$\emph{\ }onto an\emph{\ }$i$-dimensional\emph{\ }face of%
\emph{\ }$P_{2}$, for all\emph{\ }$i$, $0\leq i\leq $ dim$\left(
P_{1}\right) =$ dim$\left( P_{2}\right) $, and, in addition, $N_{P_{1}}$%
\emph{\ }onto the lattice\emph{\ }$N_{P_{2}}$, where by\emph{\ }$N_{P_{j}}$%
\emph{\ }is\emph{\ }meant the sublattice of\emph{\ }$N_{j}$
generated (as
subgroup) by aff$\left( P_{j}\right) \cap N_{j}$\emph{, }$j=1,2.$\emph{\ }If%
\emph{\ }%
\begin{equation*}
N_{1}=N_{2}=:N\ \emph{\ }\text{\ and \ rk}\left( N\right) =\emph{\ }\text{dim%
}\left( P_{1}\right) =\text{dim}\left( P_{2}\right) ,
\end{equation*}%
then these $\varpi $'s are exactly the \textit{affine integral
transformations} which are composed of unimodular transformations and $N$%
-translations.
\end{definition}

\begin{definition}[Basic simplices]
A lattice simplex is said to be \textit{basic} (or
\textit{unimodular}) if its vertices constitute a part of a\emph{\
}$\mathbb{Z}$-basis of the reference lattice (or equivalently, if
its relative, normalized volume equals $1$).
\end{definition}

\noindent Let now $U_{\sigma }=U_{\sigma ,N}$ be a $d$-dimensional
affine toric variety and $U_{\sigma ^{\prime }}=U_{\sigma ^{\prime
},N_{\sigma }}$, as in \textsf{(a)}. Assuming that $U_{\sigma }$
is Gorenstein, we may pass to another analytically isomorphic
\textquotedblleft
standard\textquotedblright\ representative as follows: Denote by $\mathbb{Z}%
^{d}$ the rectangular (standard) lattice in $\mathbb{R}^{d}$ and by $(%
\mathbb{Z}^{d})^{\vee }$ its dual lattice within
$(\mathbb{R}^{d})^{\vee }=$
Hom$_{\mathbb{R}}(\mathbb{R}^{d},\mathbb{R})$. Since dim$\left(
\sigma ^{\prime }\right) =$ rank$\left( N_{\sigma }\right) =d$, or
equivalently,
since $\left( \sigma ^{\prime }\right) ^{\vee }$ is strongly convex in $%
\left( M_{\sigma }\right) _{\mathbb{R}}$, Theorem \ref{GOR-PR}
(iii) implies
\begin{equation*}
\text{Gen}\left( \sigma ^{\prime }\right) \subset
\mathbf{H}^{\left( d\right) }\text{ \ \ with\ \ \ \
}\mathbf{H}^{\left( d\right) }:=\left\{ \mathbf{x}\in \left(
N_{\sigma }\right) _{\mathbb{R}}\ \left\vert \ \left\langle
m_{\sigma ^{\prime }},\mathbf{x}\right\rangle =1\right. \right\} ,
\end{equation*}%
for a uniquely determined $m_{\sigma ^{\prime }}\in M_{\sigma }.$ Clearly, $%
\sigma ^{\prime }\cap \mathbf{H}^{\left( d\right) }$ is a $\left(
d-1\right) $-dimensional lattice polytope (w.r.t. $N_{\sigma }$).
We choose \ a specific $\mathbb{Z}$-module isomorphism
\begin{equation*}
\Upsilon :N_{\sigma }\overset{\cong }{\longrightarrow
}\mathbb{Z}^{d}
\end{equation*}%
inducing an $\mathbb{R}$-vector space isomorphism
\begin{equation*}
\Phi =\Upsilon \otimes _{\mathbb{Z}}\text{id}_{\mathbb{R}}:\left(
N_{\sigma }\right) _{\mathbb{R}}\overset{\cong }{\longrightarrow
}\mathbb{R}^{d},
\end{equation*}%
such that
\begin{equation*}
\Phi (\mathbf{H}^{\left( d\right) })=\left\{ \mathbf{x}=\left(
x_{1},\ldots ,x_{d}\right) \in \mathbb{R}^{d}\ \left\vert \
x_{1}=1\right. \right\} \mathbf{=:\bar{H}}^{\left( d\right)
}\text{.}
\end{equation*}%
Obviously, $P:=\Phi \left( \sigma ^{\prime }\cap
\mathbf{H}^{\left( d\right) }\right) \subset
\mathbf{\bar{H}}^{\left( d\right) }$ is a lattice $\left(
d-1\right) $-dimensional polytope \ \ (w.r.t. aff$\left( P\right)
\cap \mathbb{Z}^{d}$). Defining
\begin{equation*}
\tau _{P}:=\text{pos}(P)=\left\{ \kappa \ \mathbf{x}\in \mathbb{R}^{d}%
\mathbf{\ }\left\vert \ \kappa \in \mathbb{R}_{\geq 0},\
\mathbf{x}\in P\right. \right\}
\end{equation*}%
(cf. fig. \textbf{1}) we obtain easily the following:

\begin{lemma}
\label{ISOM}\emph{(i) }There exists a torus-equivariant analytic
isomorphism
\begin{equation*}
U_{\sigma ^{\prime }}\cong U_{\tau _{P}}\ (=\text{\emph{Spec}}(\mathbb{C}%
\left[ (\mathbb{Z}^{d})^{\vee }\cap \tau _{P}^{\vee }\right] ))
\end{equation*}%
mapping \emph{orb}$\left( \sigma ^{\prime }\right) $ onto
\emph{orb}$\left( \tau _{P}\right) .$ Moreover, $U_{\tau _{P}}$ is
singular \emph{(}and its singular locus contains at least
\emph{orb}$\left( \tau _{P}\right) $\emph{)} iff $P$ is not a
basic simplex w.r.t. \emph{aff}$\left( P\right) \cap
\mathbb{Z}^{d}.\medskip $\newline \emph{(ii)} If \ $Q\subset
\mathbf{\bar{H}}^{\left( d\right) }$ is another
lattice $\left( d-1\right) $-dimensional polytope \emph{(}w.r.t. $\mathbf{%
\bar{H}}^{\left( d\right) }\cap \mathbb{Z}^{d}$\emph{)}, then $P\sim Q$ $\ $%
iff there exists a torus-equivariant analytic isomorphism $U_{\tau
_{P}}\cong U_{\tau _{Q}}$ mapping $\emph{orb}\left( \tau _{P}\right) $ onto $%
\emph{orb}\left( \tau _{Q}\right) $.
\end{lemma}

\begin{figure}[h]
\begin{center}
\includegraphics[width=8cm,height=6cm]{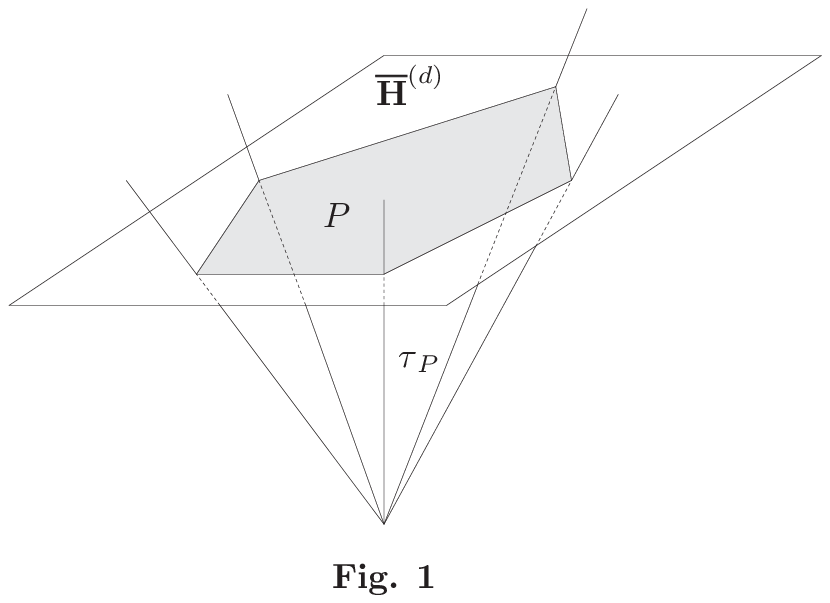}
\end{center}
\end{figure}

\noindent{}\textsf{(c) }{}Let $\mathbb{R}^{d}$ be again the usual $d$%
-dimensional euclidean space, $\mathbb{Z}^{d}$ the usual
rectangular lattice in $\mathbb{R}^{d}$ and
$(\mathbb{Z}^{d})^{\vee}$ its dual lattice. From now on, we shall
denote by $\left\{ e_{1},e_{2},\ldots,e_{d}\right\} $ the standard
$\mathbb{Z}$-basis of $\mathbb{Z}^{d}$, by $\left\{ e_{1}^{\vee
},e_{2}^{\vee},\ldots,e_{d}^{\vee}\right\} $ its dual basis, and
we shall represent the points of $\mathbb{R}^{d}$ by column
vectors and the points of its dual $(\mathbb{R}^{d})^{\vee}$ by
row vectors.

\begin{definition}
\label{FREPAR}A\textit{\ sequence\emph{\ }of free parameters of
length} $\ell $ (w.r.t.\emph{\ }$\mathbb{Z}^{d}$) is defined to be
a finite sequence
\begin{equation*}
\mathbf{m}:=\left( m_{1},m_{2},\ldots,m_{\ell}\right) ,\ \ \
1\leq\ell\leq d-1,
\end{equation*}
consisting of vectors%
\begin{equation*}
m_{i}:=\left( m_{i,1},m_{i,2},\ldots,m_{i,d}\right) ,\ 1\leq i\leq \ell,%
\text{ \ of \ }(\mathbb{Z}^{d})^{\vee}\mathbb{r}\left\{ \left(
0,\ldots,0\right) \right\}
\end{equation*}
for which\emph{\ }$m_{i,j}=0$ for all\emph{\ }$i$\emph{, }$1\leq
i\leq\ell$,
and all\emph{\ }$j$,\emph{\ }$1\leq j\leq d$, with\emph{\ }$i<j$. As\emph{\ }%
$\left( \ell\times d\right) $-matrix such an\textit{\
}$\mathbf{m}$ has the form:\smallskip\
\begin{equation}
\mathbf{m}=\left(
\begin{array}{cc}
\begin{array}{cccccc}
m_{1,1} & 0 & 0 & \cdots & \cdots & 0 \\
m_{2,1} & m_{2,2} & 0 & \cdots & \cdots & 0 \\
m_{3,1} & m_{3,2} & m_{3,3} & \cdots & \cdots & 0 \\
\vdots & \vdots & \vdots & \ddots & \cdots & \vdots \\
m_{\ell-1,1} & m_{\ell-1,2} & m_{\ell-1,3} & \cdots & \ddots & 0 \\
m_{\ell,1} & m_{\ell,2} & m_{\ell,3} & \cdots & \cdots & m_{\ell,\ell}%
\end{array}
\!\!\smallskip & \underset{d-\ell\emph{\ }\text{zero-columns}}{\underbrace{%
\begin{array}{ccc}
0 & \cdots & 0\smallskip \\
0 & \cdots & 0\smallskip \\
0 & \cdots & 0\smallskip \\
\vdots & \cdots & \vdots \\
0 & \cdots & 0 \\
0 & \cdots & 0%
\end{array}
}}%
\end{array}
\right)   \label{MATRIX-F}
\end{equation}
\end{definition}

\begin{definition}[Nakajima polytopes]
\label{NAKPOL}Fixing the dimension $d$\ of our reference space, we
define the polytopes
\begin{equation*}
\left\{ P_{\mathbf{m}}^{\left( i\right)
}\subset\mathbf{\bar{H}}^{\left( d\right) }\mathbf{\ }\left\vert \
i\in\mathbb{N}\text{\emph{, }}1\leq i\leq d\right. \right\}
\end{equation*}
lying on
\begin{equation*}
\mathbf{\bar{H}}^{\left( d\right) }\mathbf{=}\left\{ \mathbf{x}%
=(x_{1},..,x_{d})^{\intercal}\in\mathbb{R}^{d}\ \left\vert \
x_{1}=1\right. \right\}
\end{equation*}
and being associated to an \textquotedblleft
admissible\textquotedblright\
free-parameter-sequence (or matrix) $\mathbf{m}$ as in (\ref{MATRIX-F}) w.r.t%
$\emph{.\ }\mathbb{Z}^{d}$\emph{\ }(with length\emph{\ }$\ell=i-1$, for\emph{%
\ }$2\leq i\leq d$) by using induction on $i$; namely we define
\begin{equation*}
P_{\mathbf{m}}^{\left( 1\right) }:=\{e_{1}\}=\{(1,\underset{\left(
d-1\right)
\text{-times}}{\underbrace{0,0,\ldots,0,0}})^{\intercal}\},
\end{equation*}
and\emph{\ }for\emph{\ }$2\leq i\leq d$,{\small
\begin{equation}
P_{\mathbf{m}}^{\left( i\right) }:=\text{conv}\left( \left\{ P_{\mathbf{m}%
}^{\left( i-1\right) }\cup\left. \{(\mathbf{x}^{\prime
},\left\langle
m_{i-1},\mathbf{x}\right\rangle ,\underset{\left( d-i\right) \text{-times}}{%
\underbrace{0,..,0}})^{\intercal}\ \right\vert \ \mathbf{x}=(\mathbf{x}%
^{\prime},\underset{\left( d-i+1\right) \text{-times}}{\underbrace{0,..,0}%
)^{\intercal}}\in P_{\mathbf{m}}^{\left( i-1\right) }\right\}
\right)
\smallskip   \label{PIM}
\end{equation}
} where $\mathbf{x}^{\prime}=(x_{1},x_{2},..,x_{i-1}).$ $P_{\mathbf{m}%
}^{\left( i\right) }$ is obviously $\left( i-1\right) $-dimensional. For $%
\mathbf{m}$ to be
\textquotedblleft\textit{admissible}\textquotedblright\ means that
\begin{equation}
\left\langle m_{i-1},\mathbf{x}\right\rangle \geq0,\ \ \forall
\mathbf{x,\ \
\ x}=(x_{1},x_{2},\ldots,x_{i-1},\underset{\left( d-i+1\right) \text{-times}}%
{\underbrace{0,\ldots,0})^{\intercal}}\in P_{\mathbf{m}}^{\left(
i-1\right) }.   \label{ADMISSIBLE}
\end{equation}
Any lattice $\left( i-1\right) $-polytope $P$ which is lattice
equivalent to a $P_{\mathbf{m}}^{\left( i\right) }$ (as defined
above) will be called a \textit{Nakajima polytope} (w.r.t.$\emph{\
}\mathbb{R}^{d}$)\emph{. }As it is explained in \cite{DHaZ},
$P_{\mathbf{m}}^{\left( i\right) }$ can be written (for\emph{\
}$2\leq i\leq d$) as:
\begin{equation}
\left\{ \left. \mathbf{x}=(\mathbf{x}^{\prime},\underset{\left(
d-i\right)
\text{-times}}{x_{i},\underbrace{0,..,0})^{\intercal}}\in(P_{\mathbf{m}%
}^{\left( i-1\right) }\times\mathbb{R}\times\{\mathbf{0}\})\mathbb{%
\hookrightarrow}\mathbf{\bar{H}}^{\left( d\right) }\,\right\vert \
0\leq x_{i}\medskip\leq\left\langle
m_{i-1},\mathbf{x}^{\prime}\right\rangle \right\} \smallskip
\label{PRISM}
\end{equation}
i.e., as a polytope determined by means of a suitably cutted
\textquotedblleft half-line prism\textquotedblright\ over $P_{\mathbf{m}%
}^{\left( i-1\right) }.$
\end{definition}

\begin{examples}
\label{Ex}(i) For\emph{\ }$i=d=1,$\ we have trivially\emph{\ }$P_{\mathbf{m}%
}^{\left( 1\right) }=\left\{ 1\right\} .\smallskip $\emph{\newline
}(ii) For\emph{\ }$i=d=2,\ \mathbf{m}=\left( m_{1,1},0\right)
$\emph{\ }we have\emph{\ }$P_{\mathbf{m}}^{\left( 1\right)
}=\left\{ \left( 1,0\right) ^{\intercal }\right\} $\emph{\ }and
\begin{align*}
P_{\mathbf{m}}^{\left( 2\right) }& =\text{\emph{\ }conv}\left(
\left\{ \left( 1,0\right) ^{\intercal }\right\} \cup \left\{
\left( 1,\left\langle
m_{1},(1,0)\right\rangle \right) ^{\intercal }\right\} \right)  \\
& =\text{conv}\left( \left\{ \left( 1,0\right) ^{\intercal
}\right\} \cup \left\{ \left( 1,m_{1,1}\right) ^{\intercal
}\right\} \right) \text{\emph{, \ }}m_{1,1}>0.
\end{align*}%
(iii) \ For $i=d=3,$ and
\begin{equation*}
\mathbf{m=}\left(
\begin{array}{ccc}
m_{1,1} & 0 & 0 \\
m_{2,1} & m_{2,2} & 0%
\end{array}%
\right)
\end{equation*}%
we obtain
\begin{equation*}
P_{\mathbf{m}}^{\left( 3\right) }=\text{\emph{\ }conv}\left(
\left\{ \left( 1,0,0\right) ^{\intercal },\left(
1,m_{1,1},0\right) ^{\intercal },\left( 1,0,m_{2,1}\right)
^{\intercal },\left( 1,m_{1,1},m_{2,1}+m_{1,1}m_{2,2}\right)
^{\intercal }\right\} \right)
\end{equation*}%
with
\begin{equation}
m_{1,1}>0,\ \ m_{2,1}\geq 0,\ \ m_{2,1}+m_{1,1}m_{2,2}\geq 0,\ \
\left( m_{2,1},m_{2,2}\right) \neq \left( 0,0\right) .
\label{ADM3}
\end{equation}%
\newline
(iv) Finally, for\emph{\ }$i=d=4,$ and
\begin{equation*}
\mathbf{m=}\left(
\begin{array}{cccc}
m_{1,1} & 0 & 0 & 0 \\
m_{2,1} & m_{2,2} & 0 & 0 \\
m_{3,1} & m_{3,2} & m_{3,3} & 0%
\end{array}%
\right)
\end{equation*}%
we get
\begin{equation*}
P_{\mathbf{m}}^{\left( 4\right) }=\text{conv}\left( \left\{
\begin{array}{l}
\left( 1,0,0,0\right) ^{\intercal },\left( 1,m_{1,1},0,0\right)
^{\intercal
},\left( 1,0,m_{2,1},0\right) ^{\intercal }, \\
\left( 1,m_{1,1},m_{2,1}+m_{1,1}m_{2,2},0\right) ^{\intercal
},\smallskip
\left( 1,0,0,m_{3,1}\right) ^{\intercal }, \\
\left( 1,m_{1,1},0,m_{3,1}+m_{1,1}m_{3,2}\right) ^{\intercal }, \\
\left( 1,0,m_{2,1},m_{3,1}+m_{3,3}m_{2,1}\right) ^{\intercal
},\smallskip
\\
(1,m_{1,1},m_{2,1}+m_{1,1}m_{2,2}, \\
m_{3,1}+m_{3,2}m_{1,1}+m_{2,1}m_{3,3}+m_{1,1}m_{2,2}m_{3,3})^{\intercal }%
\end{array}%
\right\} \right)
\end{equation*}%
with
\begin{equation}
\left\{
\begin{array}{l}
m_{1,1}>0,\ m_{2,1}\geq 0,\ m_{2,1}+m_{1,1}m_{2,2}\geq 0,\smallskip \  \\
m_{3,1}\geq 0,\ m_{3,1}+m_{1,1}m_{3,2}\geq 0,\smallskip  \\
m_{3,1}+m_{3,3}m_{2,1}\geq
0,m_{3,1}+m_{2,1}m_{3,3}+m_{1,1}(m_{3,2}+m_{2,2}m_{3,3})\geq
0,\smallskip
\\
\left( m_{2,1},m_{2,2}\right) \neq \left( 0,0\right) ,\ \ \left(
m_{3,1},m_{3,2},m_{3,3}\right) \neq \left( 0,0,0\right) .%
\end{array}%
\right.   \label{ADM4}
\end{equation}%
In the figures \textbf{2} and \textbf{3} \ we illustrate the
lattice
polytopes\emph{\ }$P_{\mathbf{m}}^{\left( 3\right) }$, $P_{\mathbf{m}%
}^{\left( 4\right) }$, respectively, for
\begin{equation*}
\mathbf{m=}\left(
\begin{array}{rrr}
2 & 0 & 0 \\
2 & 1 & 0%
\end{array}%
\right) \text{ \ \ \ \ and \ \ \ \ }\mathbf{m=}\left(
\begin{array}{cccc}
1 & 0 & 0 & 0 \\
1 & 0 & 0 & 0 \\
2 & -1 & -1 & 0%
\end{array}%
\right) .
\end{equation*}
\end{examples}

\begin{figure}[h]
\begin{center}
\includegraphics{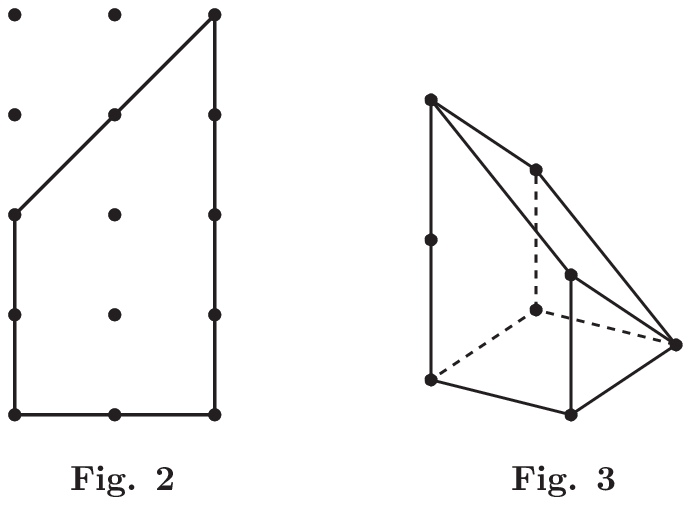}
\end{center}
\end{figure}

\begin{remark}[On the set of vertices]
A convenient reformulation of (\ref{PIM}) reads as
\begin{equation*}
P_{\mathbf{m}}^{\left( i\right) }=\text{conv}\left(
P_{\mathbf{m}}^{\left( i-1\right) }\cup
W_{i}\,P_{\mathbf{m}}^{\left( i-1\right) }\right) ,
\end{equation*}%
where
\begin{equation*}
W_{i}:=\left( e_{1}^{\vee },e_{2}^{\vee },\ldots ,e_{i-1}^{\vee
},m_{i-1},0,\ldots ,0\right) ^{\intercal },\quad i=2,\dots ,d.
\end{equation*}%
is the $(d\times d)$-matrix whose non-zero rows are $e_{1}^{\vee
},e_{2}^{\vee },\ldots ,e_{i-1}^{\vee },m_{i-1}$. Setting
\begin{equation*}
S_{\mathbf{m}}^{\left( 1\right) }:=\{e_{1}\}
\end{equation*}%
and {\small
\begin{equation*}
S_{\mathbf{m}}^{\left( i\right) }:=\left\{ \left. v+(0,\ldots
,0,\varepsilon
_{i-1}\left\langle m_{i-1},v\right\rangle ,\underset{\left( d-i\right) \text{%
-times}}{\underbrace{0,..,0}})^{\intercal }\,\right\vert
\,\varepsilon _{i}\in \{0,1\},\,v\in S_{\mathbf{m}}^{\left(
i-1\right) }\right\} ,
\end{equation*}%
}for $i=2,\ldots ,d,$ we have%
\begin{equation*}
P_{\mathbf{m}}^{\left( i\right) }=\text{ conv}\left(
S_{\mathbf{m}}^{\left( i\right) }\right) ,
\end{equation*}%
and thus, vert$(P_{\mathbf{m}}^{\left( i\right) })\subseteq S_{\mathbf{m}%
}^{\left( i\right) }$, where by $\mathrm{vert}(P)$ we denote the
set of vertices of a polytope $P$. Observe, that from the
definition of the Nakajima polytopes we know that
\begin{equation}
v\geq 0,\text{ \ for all }v\in S_{\mathbf{m}}^{\left( i\right)
},\quad i\in \{1,\dots ,d\},  \label{POSCOORD}
\end{equation}%
i.e., in particular, all the coordinates of the vertices are
non-negative.
\end{remark}

\noindent {}To provide a more explicit description of the elements
of the sets $S_{\mathbf{m}}^{\left( d\right) }$ and thus of the
possible vertices of $P_{\mathbf{m}}^{\left( d\right) }$ too, we
define for any choice of
\begin{equation*}
\varepsilon =\left( \varepsilon _{1},\ldots ,\varepsilon
_{d-1}\right) \in \{0,1\}^{d-1}
\end{equation*}%
and any $k\in \{2,\ldots ,d\}$ the vector
\begin{equation}
v_{k}(\varepsilon ):=\varepsilon _{k-1}\left(
m_{k-1,1}+\sum_{l=2}^{k-1}m_{k-1,l}v_{l}\left( \varepsilon \right)
\right) \label{RECC}
\end{equation}%
and we set
\begin{equation*}
v\left( \varepsilon \right) :=\left( 1,v_{2}\left( \varepsilon
\right) ,\ldots ,v_{d}\left( \varepsilon \right) \right)
^{\intercal }.
\end{equation*}%
On account of the definition of the sets $S_{\mathbf{m}}^{\left(
i\right) }$ we have
\begin{equation}
S_{\mathbf{m}}^{\left( d\right) }=\left\{ \left. v\left(
\varepsilon \right) \ \right\vert \ \varepsilon \in
\{0,1\}^{d-1}\right\} .  \label{ESD}
\end{equation}%
The entries $v_{2}\left( \varepsilon \right) ,\ldots ,v_{d}\left(
\varepsilon \right) $ of $v\left( \varepsilon \right) $ can be
determined by exploiting the intrinsic reccurence relation which
occurs in (\ref{RECC}). For $\varepsilon \in \{0,1\}^{d-1}$, $n\in
\{2,\ldots ,d\}$, and $k\in \{0,1,\ldots ,n-2\}$ we define the
following sum of products:

\begin{equation*}
q_{k,n}\left( \varepsilon\right)
:=\sum_{0=i_{0}<i_{1}<i_{2}<\cdots
<i_{k}<i_{k+1}=n-1}\quad\prod_{j=1}^{k+1}\varepsilon_{i_{j}}%
\,m_{i_{j},i_{j-1}+1}.
\end{equation*}
Observe that $q_{0,n}\left( \varepsilon\right) =\varepsilon_{n-1}\,m_{n-1,1}$%
.

\begin{proposition}
\label{QUS}For all $n\in\{2,\ldots,d\}$ we have
\begin{equation*}
v_{n}(\varepsilon)=\sum_{k=0}^{n-2}q_{k,n}(\varepsilon).
\end{equation*}
\end{proposition}

\noindent{}\textit{Proof}. First we notice that
\begin{equation*}
q_{k+1,n+1}(\varepsilon)=\varepsilon_{n}\sum_{i=k+2}^{n}m_{n,i}\,q_{k,i}(%
\varepsilon),
\end{equation*}
which follows immediately from the definition, because

\begin{equation}
\begin{array}{l}
\varepsilon_{n}\sum\limits_{i=k+2}^{n}m_{n,i}q_{k,i}(\varepsilon) \\
\\
=\sum\limits_{i=k+2}^{n}\varepsilon_{n}m_{n,i}\
\sum\limits_{0=j_{0}<j_{1}<j_{2}<\cdots<j_{k}<j_{k+1}=i-1}\quad\prod%
\limits_{l=1}^{k+1}\varepsilon_{j_{l}}\,m_{j_{l},j_{l-1}+1} \\
\\
=\sum\limits_{0=j_{0}<j_{1}<j_{2}<\cdots<j_{k}<j_{k+1}<j_{k+2}=n}\quad
\prod\limits_{l=1}^{k+2}\varepsilon_{j_{l}}%
\,m_{j_{l},j_{l-1}+1}=q_{k+1,n+1}(\varepsilon).%
\end{array}
\label{BIGEQUATION}
\end{equation}
We prove the Proposition by using induction w.r.t. $n$. For $n=2$
we obtain the identity
\begin{equation*}
v_{2}(\varepsilon)=\varepsilon_{1}\,m_{1,1}=q_{0,2}(\varepsilon).
\end{equation*}
Now let $n$ be $>2$. From the definition (\ref{RECC}) and our
induction hypothesis
\begin{equation*}
\begin{split}
v_{n+1}(\varepsilon) &
=\varepsilon_{n}\,m_{n,1}+\varepsilon_{n}\sum
_{k=2}^{n}m_{n,k}\,v_{k}(\varepsilon) \\
&
=q_{0,n+1}(\varepsilon)+\varepsilon_{n}\sum_{k=2}^{n}m_{n,k}\left(
\sum_{l=0}^{k-2}q_{l,k}(\varepsilon)\right) \\
&
=q_{0,n+1}(\varepsilon)+\sum_{l=0}^{n-2}\varepsilon_{n}\,%
\sum_{k=l+2}^{n}m_{n,k}q_{l,k}(\varepsilon) \\
&
=q_{0,n+1}(\varepsilon)+\sum_{l=0}^{n-2}q_{l+1,n+1}(\varepsilon),
\end{split}
\end{equation*}
where the last identity follows from (\ref{BIGEQUATION}).
\hfill$\square$

\begin{corollary}
\label{POSITIVITY}For any choice of
$\varepsilon_{\rho}\in\{0,1\}$, $\rho \in\{1,\ldots,d-2\}$, we
have
\begin{equation*}
\sum_{k=0}^{d-2}\quad\sum_{0=i_{0}<i_{1}<i_{2}<%
\cdots<i_{k}<i_{k+1}=d-1}m_{d-1,i_{k}+1}\cdot\quad\prod_{j=1}^{k}%
\varepsilon_{i_{j}}\,m_{i_{j},i_{j-1}+1}\geq0.
\end{equation*}
\end{corollary}

\noindent{}\textit{Proof}. The left hand side is nothing but
another
representation of the integer $v_{d}(\varepsilon)$, where $%
\varepsilon=\left(
\varepsilon_{1},\ldots,\varepsilon_{d-2},1\right) $,
i.e., in the case in which $\varepsilon_{d-1}=1$ (cf.~Proposition \ref{QUS}%
). In view of (\ref{ESD}) and (\ref{POSCOORD}) we have $v_{d}(\varepsilon)%
\geq0$ and the corollary is proven. \hfill$\square\bigskip$

{}\noindent{}\noindent\noindent{}\noindent{}Now Nakajima's
Classification Theorem \cite[Thm. 1.5, p. 86]{Nakajima} can be
formulated as follows:

\begin{theorem}[Nakajima's Classification of Toric L.C.I.'s]
\label{Nak-thm}Let $N$ be a free $\mathbb{Z}$-module of rank $r\geq2$, and $%
\sigma\subset N_{\mathbb{R}}$ \ a s.c.p. cone of dimension $d$,
$2\leq d\leq r$. Moreover, let $U_{\sigma}$ denote the affine
toric variety associated to
$\sigma$, and $U_{\sigma ^{\prime}}$ as in \emph{\textsf{(b)}}. Then $%
U_{\sigma}$ is a local complete intersection\emph{\ }\textbf{if and only if }%
\ there exists an admissible sequence $\mathbf{m}$ of free
parameters of length $d-1$ $\emph{(}$w.r.t.$\emph{\
}\mathbb{Z}^{d}\emph{)}$, such that\smallskip\ for any standard
representative $U_{\tau_{P}}\cong U_{\sigma^{\prime}}$ of
$U_{\sigma }$ we have $P\sim P_{\mathbf{m}}^{\left( d\right) }$,
i.e., $P$ is a Nakajima $\left( d-1\right) $-dimensional polytope
\emph{(}w.r.t. $\mathbb{R}^{d}$\emph{).}
\end{theorem}

\begin{remark}
\label{REMNAK}(i) Theorem \ref{Nak-thm} was first proved in
dimensions $2$ and $3$\ by Ishida \cite[Thm. 8.1]{Ishida}.
Previous classification results,
due to Watanabe \cite{Watanabe}, cover essentially only the class of the $%
\mathbb{Q}$-factorial toric l.c.i.'s in all dimensions. In fact,
the term
\textquotedblleft Watanabe simplex\textquotedblright\ introduced in \cite[%
5.13]{DHZ} can be used, \textit{up to lattice equivalence}, as a
synonym for a Nakajima polytope (in the sense of \ \ref{NAKPOL})
which happens to be a simplex.\medskip \newline (ii) Obviously,
$U_{\sigma }$ is a l.c.i. $\Longleftrightarrow $ $U_{\sigma
^{\prime }}\cong U_{\tau _{P}}$ is a \textquotedblleft
g.c.i\textquotedblright , i.e., a \textit{global} complete
intersection in the sense of \cite{Ishida}.\emph{\ }(It is worth
mentioning that in the setting of \ \cite{DHZ}, it was always
assumed that $d=r$; therefore, the abelian quotient
\textquotedblleft g.c.i.\textquotedblright -spaces were
abbreviated therein simply as \textquotedblleft
c.i.'s\textquotedblright ).\medskip \newline (iii) For a non-basic
Nakajima polytope $P$, $(U_{\tau _{P}}$, orb$\left( \tau
_{P}\right) )$ is a toric g.c.i.-singularity.\medskip \newline
(iv) If $P$\ is$\ $a Nakajima $\left( d-1\right) $-polytope and
$\tau _{P}$ non-basic w.r.t.\emph{\ }$\mathbb{Z}^{d}$, then the
orbit orb$\left( \tau
_{P}\right) \in U_{\tau _{P}}$ has splitting codimension $\varkappa $, with%
\emph{\ } $2\leq \varkappa \leq d-1$\emph{\ }iff\emph{\ }$P$ is
lattice-equivalent to the join\emph{\ }$\check{P}\star \mathbf{s}$ of a\emph{%
\ }$\left( \varkappa -1\right) $-dimensional (non-basic)\emph{\
}Nakajima
polytope\emph{\ }$\check{P}$ with a basic $\left( d-\varkappa -1\right) $%
-simplex\emph{\ }$\mathbf{s}$.\medskip \emph{\newline }(v) It is
easy for every\emph{\ }$P\subset \mathbf{\bar{H}}^{\left( d\right)
}$, with\emph{\ }$P\sim P_{\mathbf{m}}^{\left( d\right) }$, to
verify that
\begin{equation}
d\leq \#(\text{vert}(P))\leq 2^{d-1}\text{ \ \ and \ }d\leq \#(\{\text{%
facets of\emph{\ }}P\})\leq 2\left( d-1\right) .  \label{BOUNDS}
\end{equation}%
(vi) The question:
\begin{equation*}
\text{\textquotedblleft what kind of equations define toric
l.c.i.-singularities \ }(U_{\tau _{P}},\text{orb}\left( \tau
_{P}\right) )?\textquotedblright
\end{equation*}
will be answered only in dimensions $\geq 3$, because, as it is
well-known, in dimension $2$ only the classical \textquotedblleft
Kleinian\textquotedblright\ hypersurface singularities
\begin{equation*}
\left\{ \left. \left( z,w,t\right) \in \mathbb{C}^{3}\ \right\vert
\ z^{\kappa }-wt=0\right\} ,\ \ \kappa \in \mathbb{Z}_{\geq 2},
\end{equation*}%
of \textquotedblleft type $\mathbf{A}_{\kappa -1}"$ are present.
\end{remark}\bigskip

\section{Equations Defining Toric L.C.I.-Singularities}

\noindent{}Our main result is the following: Let $N$ be a free $\mathbb{Z}$%
-module of rank $r\geq3$, and $\sigma\subset N_{\mathbb{R}}$ \ a
rational s.c.p.c. of dimension $d$, $3\leq d\leq r$, such that
$U_{\sigma}=U_{\sigma ,N}$ is a local complete
intersection.\emph{\ }By Theorem \ref{Nak-thm} there exists an
admissible sequence
\begin{equation}
\mathbf{m}=\left(
\begin{array}{ccccccc}
m_{1,1} & 0 & 0 & \cdots & \cdots & 0 & 0 \\
m_{2,1} & m_{2,2} & 0 & \cdots & \cdots & 0 & 0 \\
m_{3,1} & m_{3,2} & \ddots & \cdots & \cdots & \vdots & \vdots \\
\vdots & \vdots & \ddots & \ddots &  & \vdots & \vdots \\
\vdots & \vdots & \vdots & \ddots & \ddots & 0 & 0 \\
m_{d-1,1} & m_{d-1,2} & m_{d-1,3} & \cdots & m_{d-1,d-2} & m_{d-1,d-1} & 0%
\end{array}
\!\!\right)   \label{NAKMATR}
\end{equation}
of free parameters of length $d-1$ (w.r.t. the lattice$\emph{\ }\mathbb{Z}%
^{d}$), such that\smallskip\ for any standard representative $%
U_{\tau_{P}}\cong U_{\sigma^{\prime}}$ of $U_{\sigma}$ we have $P\sim P_{%
\mathbf{m}}^{\left( d\right) }$\emph{. }We may, without loss of
generality, assume that $(U_{\tau_{P}}$, orb$\left(
\tau_{P}\right) )$ is an msc-singularity.

\begin{theorem}
\label{MAIN}If $P\sim P_{\mathbf{m}}^{\left( d\right) },$ as
above, then
\begin{equation*}
U_{\tau_{P}}\cong U_{\tau_{P_{\mathbf{m}}^{\left( d\right) }}}\cong\text{%
\emph{Spec}}\left( \mathbb{C}\left[ z_{1},z_{2},\ldots
,z_{d},z_{d+1},\ldots,z_{2d-1}\right] \ /\ \mathcal{I}\right)
\end{equation*}
with defining ideal\smallskip\
\begin{equation*}
\fbox{$%
\begin{array}{ccc}
&  &  \\
& \mathcal{I}=\left( \left\{ \left.
{\displaystyle\prod\limits_{1\leq i\leq
j}}z_{i}{}^{\mu_{i-1,j}\smallskip\smallskip}\,z_{d+i-1}^{-\lambda
_{i-1,j}\smallskip}{}\,-\,z_{j+1}z_{d+j}\ \right| \ 1\leq j\leq
d-1\right\}
\right) &  \\
&  &
\end{array}
$}\smallskip
\end{equation*}
where $\lambda_{i,j}$'s, $j\in\{1,\ldots,d-1\}$, are non-positive
integers determined recursively by the formula\smallskip\
\begin{equation}
\lambda_{i,j}:=\left\{
\begin{array}{ll}
\min\left\{ 0,m_{j,i+1}\right\} ,\  & \text{if }i=j-1,\smallskip \\
\min\left\{
0,m_{j,i+1}+\sum\limits_{k=i+1}^{j-1}m_{k,i+1}\,\lambda
_{k,j}\right\} , & \text{if }i\in\{j-2,\ldots,1\},\smallskip \\
0, & \text{if }i=0,%
\end{array}
\right. \smallskip   \label{LAMDAS}
\end{equation}
and $\mu_{i,j}$'s, $j\in\{1,\ldots,d-1\}$, are non-negative
integers defined
by the formula\smallskip%
\begin{equation}
\mu_{i,j}:=\left\{
\begin{array}{ll}
\max\left\{ 0,m_{j,i+1}\right\} ,\  & \text{if }i=j-1,\smallskip \\
\max\left\{
0,m_{j,i+1}+\sum\limits_{k=i+1}^{j-1}m_{k,i+1}\,\lambda
_{k,j}\right\} , & \text{if }i\in\{j-2,\ldots,1\},\smallskip \\
m_{j,1}+\sum\limits_{k=1}^{j-1}m_{k,1}\lambda_{k,j}, & \text{if }i=0.%
\end{array}
\right.   \label{MIOUS}
\end{equation}
\end{theorem}

\begin{remark}
\label{REMMAIN}(i) In the formula (\ref{LAMDAS}) the
$\lambda_{j-1,j}$'s are
known from the beginning. For all $\rho\in\{2,\ldots,j-1\}$, the $%
\lambda_{j-\rho,j}$'s are to be found successively by means of
integer
linear combinations of $\lambda_{j-1,j}$, $\lambda_{j-2,j}$,$\ldots$, $%
\lambda _{j-\rho-1,j}$ (with known
coefficients).\smallskip\newline
(ii) For all $j\in\{1,\ldots,d-1\}$ and $i\in\{1,\ldots,j\}$, either $%
\mu_{i-1,j}=0$ or $\lambda_{i-1,j}=0$ (by definition). Hence, the
first monomial of each of the $d-1$ binomials which generate
$\mathcal{I}$
contains only one of the two variables $z_{i}$, $z_{d+i-1}$, $1\leq i\leq j$%
.\smallskip\newline (iii) If \textit{all} entries in
(\ref{NAKMATR}) are non-negative, then
\begin{equation}
\mathcal{I}=\left( \left\{ \left. \prod_{1\leq i\leq
j}z_{i}^{m_{j,i}}-z_{j+1}\,z_{d+j}\,\ \right| \ 1\leq j\leq
d-1\right\} \right) ,   \label{SPECIALCASE}
\end{equation}
because in this case all exponents $\lambda_{i-1,j}$ are $=0$ and
$\mu _{i-1,j}=m_{j,i}.$
\end{remark}

{}\noindent{}Next, we define
\begin{equation}
\mathfrak{Q}_{\mathbf{m}}^{\left( d\right) }:=\left\{ k\in\left\{
1,..,d-1\right\} \ \left\vert
\begin{array}{l}
\ e_{k}^{\vee}=m_{\gamma_{k}},\text{ } \\
\text{for some\smallskip\ index } \\
\gamma_{k}\in\{k,..,d-1\}%
\end{array}
\right. \right\}   \label{QINDEX}
\end{equation}
and
\begin{equation}
\mathfrak{R}_{\mathbf{m}}^{\left( d\right) }:=\left\{ l\in\left\{
1,..,d-2\right\} \ \left\vert
\begin{array}{l}
\ m_{l}-e_{l+1}^{\vee}=m_{\delta_{l}},\text{ } \\
\text{for some\smallskip\ index } \\
\delta_{l}\in\{l+1,..,d-1\}%
\end{array}
\right. \right\} .   \label{RINDEX}
\end{equation}

\begin{corollary}
\label{KOROLLAR}If $P\sim P_{\mathbf{m}}^{\left( d\right) },$ as in \emph{%
Thm. \ref{MAIN}}, then $U_{\tau_{P}}\cong
U_{\tau_{P_{\mathbf{m}}^{\left( d\right) }}}$ admits the
\emph{\textquotedblleft minimal\textquotedblright} embedding
\begin{equation*}
U_{\tau_{P_{\mathbf{m}}^{\left( d\right) }}}\hookrightarrow\mathbb{C}^{\#(%
\mathbf{Hilb}_{(\mathbb{Z}^{d})^{\vee}}(\tau_{P_{\mathbf{m}}^{\left(
d\right) }}^{\vee}))}
\end{equation*}
after eliminating the redundant variables of $\mathbb{C}\left[
z_{1},z_{2},\ldots,z_{d},z_{d+1},\ldots,z_{2d-1}\right] $. More precisely,%
{\small \smallskip%
\begin{equation*}
U_{\tau_{P_{\mathbf{m}}^{\left( d\right)
}}}\cong\text{\emph{Spec}}\left( \mathbb{C}\left[ \left\{
z_{k}\,\left\vert \,k\in\left\{ 1,..,d\right\}
\mathbb{r}\mathfrak{Q}_{\mathbf{m}}^{\left( d\right) }\right.
\right\}
\cup\left\{ z_{d+l}\,\left\vert \,l\in\left\{ 1,..,d-1\right\} \mathbb{r}%
\mathfrak{R}_{\mathbf{m}}^{\left( d\right) }\right. \right\}
\right] \ /\ \mathcal{I}\right) ,
\end{equation*}
}\smallskip where the ideal $\mathcal{I}$ \ is generated by $\#(\mathbf{Hilb}%
_{(\mathbb{Z}^{d})^{\vee}}(\tau_{P_{\mathbf{m}}^{\left( d\right)
}}^{\vee }))-d$ binomials\emph{. }These binomials are exactly
those remaining from
\begin{equation*}
\left\{ \left. {\displaystyle\prod\limits_{1\leq i\leq j}}%
z_{i}{}^{\mu_{i-1,j}\smallskip\smallskip}\,z_{d+i-1}^{-\lambda_{i-1,j}%
\smallskip}{}\,-\,z_{j+1}z_{d+j}\ \right\vert \ 1\leq j\leq d-1\right\}
\end{equation*}
after the elimination of \ the variables
\begin{equation*}
\left\{ z_{k}\ \left\vert \ k\in\mathfrak{Q}_{\mathbf{m}}^{\left(
d\right)
}\right. \right\} \cup\left\{ z_{d+l}\ \left\vert \ l\in\mathfrak{R}_{%
\mathbf{m}}^{\left( d\right) }\right. \right\}
\end{equation*}
by means of the substitutions
\begin{equation*}
z_{k}=z_{\gamma_{k}+1}\,z_{d+\gamma_{k}}\text{ \ \ \ and\ \ \ \ }%
z_{d+l}=z_{\delta_{l}+1}\,z_{d+\delta_{l}}\ .
\end{equation*}
\end{corollary}

\noindent{}Proofs of \ref{MAIN} and \ref{KOROLLAR} are given in
the next section. Let us first apply them to a couple of examples
of Nakajima polytopes.

\begin{examples}
\label{BEISPIELE}(i) For the Nakajima quadrilateral
$P_{\mathbf{m}}^{\left( 3\right) }$ of figure \textbf{2, }
(\ref{SPECIALCASE}) gives:
\begin{equation*}
U_{\tau _{P_{\mathbf{m}}^{\left( 3\right) }}}\cong \text{Spec}\left( \mathbb{%
C}\left[ z_{1},z_{2},z_{3},z_{4},z_{5}\right] \ /\ \left(
z_{1}^{2}-z_{2}z_{4},\,z_{1}^{2}z_{2}-z_{3}z_{5}\right) \right)
.\smallskip
\newline
\medskip
\end{equation*}%
(ii) Let $P_{\mathbf{m}}^{\left( 3\right) }$ be the Nakajima
triangle with
\begin{equation*}
\mathbf{m}=\left(
\begin{array}{ccc}
k & 0 & 0 \\
k & -1 & 0%
\end{array}%
\right) ,\ \ \ k\in \mathbb{Z}_{\geq 2}.
\end{equation*}%
\newline
Eliminating the variable $z_{4}\left( =z_{3}z_{5}\right) $ as in Cor. \ref%
{KOROLLAR}, and setting $w=z_{1}$, $t_{1}=z_{2}$, $t_{2}=z_{3}$,
$t_{3}=z_{5} $, we obtain the hypersurface
\begin{equation*}
U_{\tau _{P_{\mathbf{m}}^{\left( 3\right) }}}\cong \text{Spec}\left( \mathbb{%
C}\left[ w,t_{1},t_{2},t_{3}\right] \ /\ \left(
w^{k}-t_{1}t_{2}t_{3}\right) \right) .\smallskip \newline
\medskip
\end{equation*}%
(iii) For the Nakajima solid $P_{\mathbf{m}}^{\left( 4\right) }$
of figure
\textbf{3}, we have%
\begin{equation*}
\left\{
\begin{array}{lll}
\mu _{0,1}=1, & \mu _{0,2}=1, & \mu _{0,3}=0, \\
\mu _{1,2}=0, & \mu _{1,3}=0, & \mu _{2,3}=0, \\
\lambda _{0,1}=0, & \lambda _{0,2}=0, & \lambda _{0,3}=0, \\
\lambda _{1,2}=0, & \lambda _{1,3}=-1, & \lambda _{2,3}=-1.%
\end{array}%
\right.
\end{equation*}%
Hence, Theorem \ref{MAIN} gives%
\begin{equation*}
U_{\tau _{P_{\mathbf{m}}^{\left( 4\right) }}}\cong \text{Spec}\left( \mathbb{%
C}\left[ z_{1},z_{2},z_{3},z_{4},z_{5},z_{6},z_{7}\right] \ /\
\left(
z_{1}-z_{2}z_{5},z_{1}-z_{3}z_{6},z_{5}z_{6}-z_{4}z_{7}\right)
\right)
\end{equation*}%
Therefore it is possible (by Cor. \ref{KOROLLAR}) to erase the
redundant
variable $z_{1}$ (cf. (\ref{QINDEX})) and describe $U_{\tau _{P_{\mathbf{m}%
}^{\left( 4\right) }}}$ as complete intersection of two binomials in $%
\mathbb{C}_{\left( t_{1},\ldots ,t_{6}\right) }^{6}$\smallskip , where $%
t_{i}=z_{i+1}$, $1\leq i\leq 6$, as follows:
\begin{equation*}
U_{\tau _{P_{\mathbf{m}}^{\left( 4\right) }}}\cong \text{Spec}\left( \mathbb{%
C}\left[ t_{1},t_{2},t_{3},t_{4},t_{5},t_{6}\right] \ /\ \left(
t_{1}t_{4}-t_{2}t_{5},t_{4}t_{5}-t_{3}t_{6}\right) \right) .\
\end{equation*}%
(iv) Let us now give an example of a Nakajima polytope with the
smallest number of vertices (cf. (\ref{BOUNDS})), generalizing
slightly (ii). For$\ k\in \mathbb{Z}_{\geq 2}$, $d\in
\mathbb{Z}_{\geq 4}$, let
\begin{equation*}
\mathbf{s}_{k}^{\left( d\right) }\subset \mathbf{\bar{H}}^{\left(
d\right) }\hookrightarrow \mathbb{R}^{d}
\end{equation*}%
denote the $\left( d-1\right) $-simplex
\begin{equation*}
\mathbf{s}_{k}^{\left( d\right) }:=\text{conv}\left( \left\{
e_{1},e_{1}+k\,e_{2},e_{1}+k(\,e_{2}+e_{3}),\ldots
,e_{1}+k(\,e_{2}+e_{3}+\cdots +e_{d})\right\} \right)
\end{equation*}%
being constructed by the $k$-th dilation of a basic $\left( d-1\right) $%
-simplex. Obviously,
\begin{equation*}
\mathbf{s}_{k}^{\left( d\right) }=P_{\mathbf{m}}^{\left( d\right)
}
\end{equation*}%
with $\mathbf{m}$ denoting the $\left( \left( d-1\right) \times d\right) $%
-matrix having entries $m_{1,1}=k$, $m_{i,i}=1$ in its diagonal, $\forall i$%
, $2\leq i\leq d-1$, and zero entries otherwise. By Thm.
\ref{MAIN} we can embed $U_{\tau _{\mathbf{s}_{k}^{\left( d\right)
}}}$ into $\mathbb{C}^{2d-1} $ via the $d-1$ equations
\begin{equation*}
\left\{
\begin{array}{l}
z_{1}^{k}-z_{2}z_{d+1}=0,\smallskip  \\
z_{2}-z_{3}z_{d+2}=0, \\
\vdots  \\
z_{d-1}-z_{d}z_{2d-1}=0.%
\end{array}%
\right.
\end{equation*}%
In fact, replacing succesively $z_{2}$ by $z_{3}z_{d+2},$ $z_{3}$ by $%
z_{4}z_{d+3}$ etc. in the first equation (according to the pattern of Cor. %
\ref{KOROLLAR}) and setting $w=z_{1},t_{i}=z_{d+i-1}$, for all
$i\in \{1,\ldots ,d\}$, we may represent $U_{\tau
_{\mathbf{s}_{k}^{\left( d\right) }}}$ as a hypersurface embedded
into $\mathbb{C}_{\left( w,t_{1},t_{2},t_{3},\ldots ,t_{d}\right)
}^{d+1}$:
\begin{equation*}
U_{\tau _{\mathbf{s}_{k}^{\left( d\right) }}}\cong \text{Spec}\left( \mathbb{%
C}\left[ w,t_{1},t_{2},t_{3},\ldots ,t_{d}\right] \,/\,(w^{k}-\prod%
\limits_{j=1}^{d}t_{j})\right) .
\end{equation*}%
(Notice that $U_{\tau _{\mathbf{s}_{k}^{\left( d\right) }}}\cong \mathbb{C}%
^{d}/G\left( d;k\right) $ is an abelian quotient space w.r.t. a group $%
G\left( d;k\right) \cong \left(
\mathbb{Z\,}/\,k\,\mathbb{Z}\right) ^{d-1}$,
cf. \cite[Ex. 1.5, p. 90]{Watanabe} and \cite[Prop. 5.10, p. 217]{DHZ}%
.)\medskip \newline
(v) Finally, let $k_{1},k_{2},\ldots ,k_{d-1}$ be a $\left( d-1\right) $%
-tuple of positive integers ($d\geq 4$), with $k_{1}\geq 2$, and
let
\begin{equation*}
\begin{array}{ll}
\mathbf{RP\medskip }\left( k_{1},k_{2},..,k_{d-1}\right)  &
=\left\{ \left( x_{1},..,x_{d}\right) ^{\intercal }\in
\mathbb{R}^{d}\ \left\vert
\begin{array}{c}
\ x_{1}=1\medskip ,\ 0\leq x_{j+1}\leq k_{j},\  \\
\forall j,\ 1\leq j\leq d-1%
\end{array}%
\right. \right\}  \\
& \  \\
\, & =\left\{ 1\right\} \times \left[ 0,k_{1}\right] \times \left[ 0,k_{2}%
\right] \times \cdots \times \left[ 0,k_{d-1}\right]
\end{array}%
\end{equation*}%
denote the $\left( d-1\right) $-dimensional \textit{rectangular
parallelepiped} in $\mathbf{\bar{H}}^{\left( d\right)
}\hookrightarrow
\mathbb{R}^{d}$ having them as lengths of its edges (cf. \cite{DHaZ}%
).\medskip
\begin{figure}[h]
\begin{center}
\includegraphics[width=10cm,height=7cm]{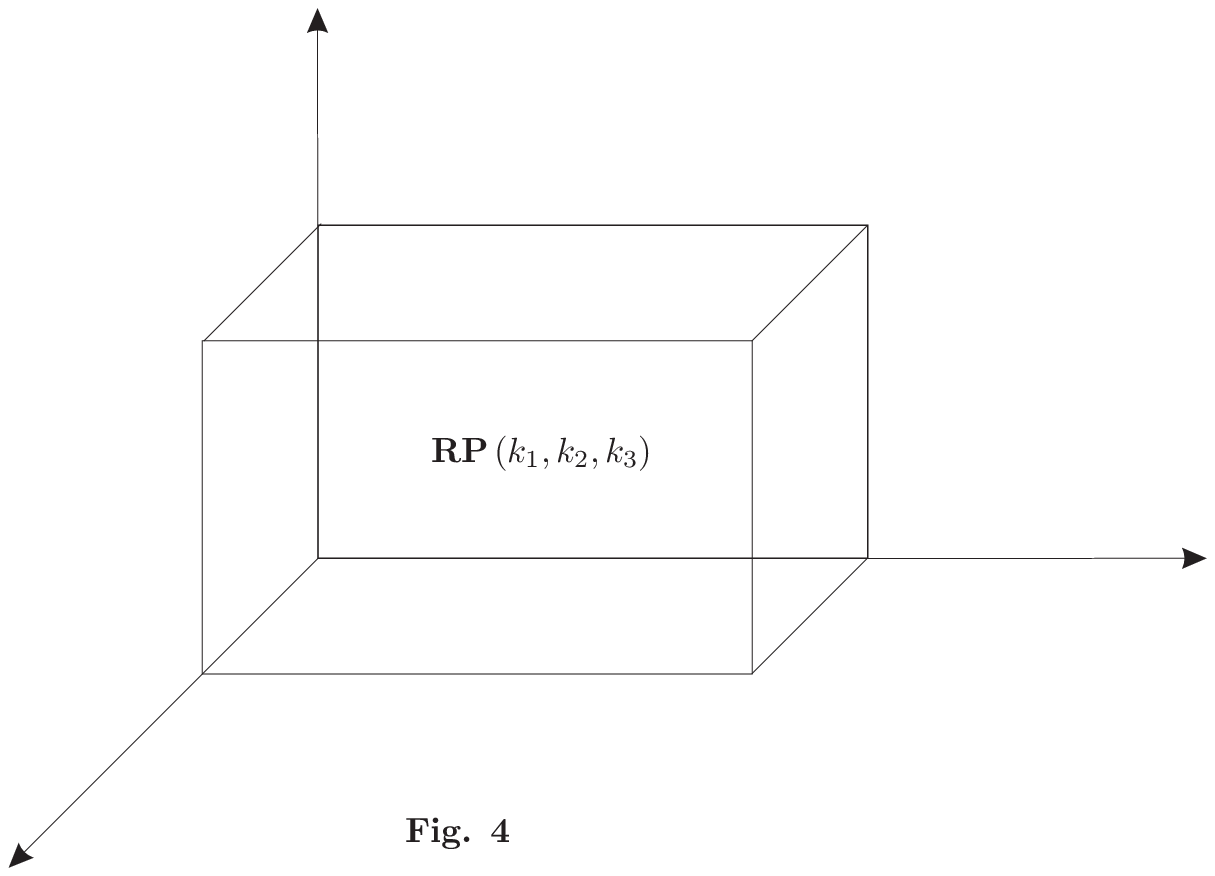}
\end{center}
\end{figure}
\newline
$\mathbf{RP\medskip }\left( k_{1},k_{2},\ldots ,k_{d-1}\right) $
has $2^{d-1}
$ vertices (i.e., the greatest possible number of vertices, cf. (\ref{BOUNDS}%
)), namely
\begin{equation*}
\{\left. e_{1}+\varepsilon _{1}\cdot k_{1}\cdot e_{2}+\varepsilon
_{2}\cdot k_{2}\cdot e_{3}+\cdots +\varepsilon _{d-1}\cdot
k_{d-1}\cdot e_{d}\ \right\vert \ \varepsilon _{1},..,\varepsilon
_{d-1}\in \{0,1\}\}
\end{equation*}%
and equals $P_{\mathbf{m}}^{\left( d\right) }$ with $\mathbf{m}$
denoting
the $((d-1)\times d)$-matrix with entries $m_{i,1}=k_{i}$, for all $i$, $%
1\leq i\leq d-1$, in its first column, and all the other entries $=0$. By (%
\ref{SPECIALCASE}) we obtain
\begin{equation*}
U_{\tau _{\mathbf{RP}\left( k_{1},..,k_{d-1}\right) }}\cong \text{Spec}(%
\mathbb{C}\left[ z_{1},..,z_{2d-1}\right] \,/\,(%
\{z_{j}^{k_{j}}-z_{j+1}z_{d+j}\ \left\vert \ 1\leq j\leq
d-1\right. \})).\medskip
\end{equation*}
\end{examples}\bigskip

\section{Proof of the Main Theorem}

\noindent{}To prove Theorem \ref{MAIN} we need several auxiliary
Lemmas. At
first, starting with an admissible sequence $\mathbf{m}$ as in (\ref{NAKMATR}%
), where $d\geq3$, we define the set
\begin{equation*}
\mathcal{L}_{\mathbf{m}}^{\left( d\right) }:=\left\{
e_{1}^{\vee}\right\} \cup\left\{ \left.
e_{k}^{\vee},m_{k-1}-e_{k}^{\vee}\ \right| \ \ 2\leq k\leq
d\right\} .
\end{equation*}

\begin{lemma}
\label{LEMMA0}$\tau_{P_{\mathbf{m}}^{\left( d\right) }}^{\vee}=$ \emph{pos}$%
\left( \mathcal{L}_{\mathbf{m}}^{\left( d\right) }\right) .$
\end{lemma}

\noindent{}\textit{Proof}. See Nakajima \cite[p. 92]{Nakajima}.\hfill $%
\square\smallskip\bigskip$

\noindent {}For $\varepsilon _{i}\in \left\{ 0,1\right\} ,$ $i\in
\left\{ 1,\ldots ,d-1\right\} ,$ we define the s.c.p. cones
$C_{\varepsilon
_{1},\ldots ,\varepsilon _{d-1}}^{\left( d\right) }\subset (\mathbb{R}%
^{d})^{\vee }$ as follows:
\begin{equation*}
C_{\varepsilon _{1},..,\varepsilon _{d-1}}^{\left( d\right) }:=\text{pos}%
\left( \{e_{1}^{\vee }\}\cup \{\left. \varepsilon
_{i}\,e_{i+1}^{\vee }+\left( 1-\varepsilon _{i}\right) \left(
m_{i}-e_{i+1}^{\vee }\right) \ \right\vert \ 1\leq i\leq
d-1\}\right) .
\end{equation*}

\begin{lemma}
\label{LEMMA1}The $2^{d-1}$ s.c.p. cones $\left\{ \left.
C_{\varepsilon _{1},\ldots,\varepsilon_{d-1}}^{\left( d\right) }\
\right| \ \varepsilon _{i}\in\left\{ 0,1\right\} ,\ 1\leq i\leq
d-1\right\} $ form a subdivision of \
$\tau_{P_{\mathbf{m}}^{\left( d\right) }}^{\vee}$ into basic cones
w.r.t. $(\mathbb{Z}^{d})^{\vee}.$
\end{lemma}

\noindent{}\textit{Proof}. Obviously, we have:
\begin{equation*}
\left\vert \text{det}\left(
e_{1}^{\vee},\varepsilon_{1}\,e_{2}^{\vee }+\left(
1-\varepsilon_{1}\right) \left( m_{1}-e_{2}^{\vee}\right)
,\ldots,\varepsilon_{d-1}\,e_{d}^{\vee}+\left(
1-\varepsilon_{d-1}\right) \left( m_{d-1}-e_{d}^{\vee}\right)
\right) \right\vert =1,
\end{equation*}
which means that all the cones $C_{\varepsilon_{1},\ldots,%
\varepsilon_{d-1}}^{\left( d\right) }$ are basic w.r.t. $(\mathbb{Z}%
^{d})^{\vee}.$ Next, we show that the intersection of two of these
simplicial cones, say of $C_{\varepsilon_{1},\ldots,\varepsilon_{d-1}}^{%
\left( d\right) }$ and $C_{\varepsilon_{1}^{\prime},\ldots,%
\varepsilon_{d-1}^{\prime}}^{\left( d\right) }$, is either a face
of both or empty. More precisely, we shall prove that
\begin{equation}
\begin{array}{l}
C_{\varepsilon_{1},\ldots,\varepsilon_{d-1}}^{\left( d\right)
}\cap
C_{\varepsilon_{1}^{\prime},\ldots,\varepsilon_{d-1}^{\prime}}^{\left(
d\right) } \\
\  \\
=\text{pos}\left( \left\{ e_{1}^{\vee}\right\} \cup\left\{
\varepsilon _{i-1}\,e_{i}^{\vee}+\left( 1-\varepsilon_{i-1}\right)
\left( m_{i-1}-e_{i}^{\vee}\right) \ \left\vert
\begin{array}{l}
\text{for all } \\
i\in\{2,..,d\}\text{ } \\
\text{for which \ } \\
\varepsilon_{i-1}=\varepsilon_{i-1}^{\prime}%
\end{array}
\right. \right\} \right) .%
\end{array}
\label{INTERSE}
\end{equation}
The inclusion \textquotedblleft$\supseteq$\textquotedblright\ is
obvious. For every element $c\in
C_{\varepsilon_{1},\ldots,\varepsilon_{d-1}}^{\left( d\right)
}\cap
C_{\varepsilon_{1}^{\prime},\ldots,\varepsilon_{d-1}^{\prime}}^{\left(
d\right) }$ there exist
\begin{equation*}
\left( \nu_{1},\ldots,\nu_{d}\right) ,\ \left( \xi_{1},\ldots,\xi
_{d}\right) \in\left( \mathbb{R}_{\geq0}\right) ^{d},
\end{equation*}
such that
\begin{equation}
\begin{array}{r}
c=\nu_{1}e_{1}^{\vee}+\sum\limits_{i=2}^{d}\nu_{i}\left(
\varepsilon _{i-1}\,e_{i}^{\vee}+\left( 1-\varepsilon_{i-1}\right)
\left(
m_{i-1}-e_{i}^{\vee}\right) \right) \\
\  \\
=\xi_{1}e_{1}^{\vee}+\sum\limits_{i=2}^{d}\xi_{i}\left(
\varepsilon _{i-1}^{\prime}\,e_{i}^{\vee}+\left(
1-\varepsilon_{i-1}^{\prime}\right)
\left( m_{i-1}-e_{i}^{\vee}\right) \right) .%
\end{array}
\label{CSXE}
\end{equation}
For the last coordinate of $c$ we obtain
\begin{equation*}
\nu_{d}\left( \varepsilon_{d-1}\,e_{d}^{\vee}+\left( 1-\varepsilon
_{d-1}\right) \left( -e_{d}^{\vee}\right) \right) =\xi_{d}\left(
\varepsilon_{d-1}^{\prime}\,e_{d}^{\vee}+\left(
1-\varepsilon_{d-1}^{\prime }\right) \left( -e_{d}^{\vee}\right)
\right) ,
\end{equation*}
and thus either $\nu_{d}=\xi_{d}=0$ or $\nu_{d}=\xi_{d}>0$ and
$\varepsilon _{d-1}=\varepsilon_{d-1}^{\prime}.$ Hence,
(\ref{CSXE}) is reduced to
\begin{equation*}
\begin{array}{r}
\nu_{1}e_{1}^{\vee}+\sum\limits_{i=2}^{d-1}\nu_{i}\left(
\varepsilon _{i-1}\,e_{i}^{\vee}+\left( 1-\varepsilon_{i-1}\right)
\left(
m_{i-1}-e_{i}^{\vee}\right) \right) \\
\  \\
=\xi_{1}e_{1}^{\vee}+\sum\limits_{i=2}^{d-1}\xi_{i}\left(
\varepsilon _{i-1}^{\prime}\,e_{i}^{\vee}+\left(
1-\varepsilon_{i-1}^{\prime}\right)
\left( m_{i-1}-e_{i}^{\vee}\right) \right) .%
\end{array}
\end{equation*}
Applying the same argumentation to the $(d-1)$-th coordinate of
this vector and repeating (actually by induction w.r.t. $d$, cf.
(\ref{PRISM})) the whole procedure, we get
\begin{equation*}
(\nu_{i}=\xi_{i}=0)\text{ \ or \ \ }(\nu_{i}=\xi_{i}>0\ \ \text{and \ }%
\varepsilon_{i-1}=\varepsilon_{i-1}^{\prime})\,,\ \text{for
}i\in\left\{ 2,\ldots,d\right\} ,
\end{equation*}
which implies the \textquotedblleft$\subseteq$\textquotedblright-part of (%
\ref{INTERSE}). Finally, we use induction w.r.t. $d$ to prove the
equality
\begin{equation}
\bigcup\limits_{\varepsilon_{i}\in\left\{ 0,1\right\} ,\ 1\leq
i\leq
d-1}C_{\varepsilon_{1},\ldots,\varepsilon_{d-1}}^{\left( d\right) }=\tau_{P_{%
\mathbf{m}}^{\left( d\right) }}^{\vee}.   \label{TCE}
\end{equation}
For $d=3$ this was proved in \cite[Lemma 8.12, p. 142]{Ishida}.
Suppose that
$d$ is $>3.$ Clearly, the inclusion \textquotedblleft$\subseteq$%
\textquotedblright\ in (\ref{TCE}) is always valid (by Lemma
\ref{LEMMA0}). Now let $v$ be an arbitrary element of
$\tau_{P_{\mathbf{m}}^{\left(
d\right) }}^{\vee}$ and $\nu_{1},\ldots,\nu_{d}$, $\xi_{1},\ldots,\xi_{d}\in%
\mathbb{R}_{\geq0},$ such that
\begin{equation*}
v=\left(
\nu_{1}e_{1}^{\vee}+\sum\limits_{i=2}^{d-1}\nu_{i}e_{i}^{\vee}+\xi_{i}\left(
m_{i-1}-e_{i}^{\vee}\right) \right)
+\nu_{d}e_{d}^{\vee}+\xi_{d}\left( m_{d-1}-e_{d}^{\vee}\right) .
\end{equation*}
Define
\begin{equation*}
\widetilde{v}:=\nu_{1}e_{1}^{\vee}+\sum\limits_{i=2}^{d-1}\nu_{i}e_{i}^{\vee
}+\xi_{i}\left( m_{i-1}-e_{i}^{\vee}\right) .
\end{equation*}
Then $\widetilde{v}\in\tau_{P_{\widetilde{\mathbf{m}}}^{\left(
d-1\right) }}^{\vee}\subset\tau_{P_{\mathbf{m}}^{\left( d\right)
}}^{\vee},$ with
\begin{equation*}
\mathbf{m}=\left(
\begin{array}{cc}
\widetilde{\mathbf{m}}\ \ 0 & 0 \\
m_{d-1} & 0%
\end{array}
\right)
\end{equation*}
and we may write
\begin{equation*}
v=\left\{
\begin{array}{ll}
\widetilde{v}+\xi_{d}m_{d-1}+\left( \nu_{d}-\xi_{d}\right)
e_{d}^{\vee}, &
\text{if \ }\nu_{d}\geq\xi_{d}, \\
\  &  \\
\widetilde{v}+\nu_{d}m_{d-1}+\left( \xi_{d}-\nu_{d}\right) \left(
m_{d-1}-e_{d}^{\vee}\right) , & \text{if \ }\nu_{d}\leq\xi_{d}.%
\end{array}
\right.
\end{equation*}
Since we know that $\widetilde{\mathbf{m}}$ is admissible, we have $%
\widetilde{v}+rm_{d-1}\in$
$\tau_{P_{\widetilde{\mathbf{m}}}^{\left( d-1\right) }}^{\vee}$
for any $r\in\mathbb{R}_{\geq0}$. Thus, in the case in which
$\nu_{d}\geq\xi_{d}$, there exists, by induction hypothesis, a
$\left(
d-1\right) $-dimensional cone $C_{\varepsilon_{1},\ldots,%
\varepsilon_{d-2}}^{\left( d-1\right) }\subset\tau_{P_{\widetilde{\mathbf{m}}%
}^{\left( d-1\right) }}^{\vee}$ containing
$\widetilde{v}+\xi_{d}m_{d-1}$ . But this shows that $v\in
C_{\varepsilon_{1},\ldots,\varepsilon_{d-2},1}^{\left( d\right)
}.$ In the other case, i.e., whenever $\nu_{d}\leq\xi_{d}$, we
find in the same way a cone of type
$C_{\varepsilon_{1},\ldots,\varepsilon_{d-2},0}^{\left( d\right)
}$ containing $v$.\hfill$\square$

\begin{lemma}
\label{LEMMA2}The set $\mathcal{L}_{\mathbf{m}}^{\left( d\right)
}$ is a system of generators of the additive semigroup
$\tau_{P_{\mathbf{m}}^{\left( d\right)
}}^{\vee}\cap(\mathbb{Z}^{d})^{\vee}.$
\end{lemma}

\noindent{}\textit{Proof}. For the subdivision of $\tau_{P_{\mathbf{m}%
}^{\left( d\right) }}^{\vee}$ into \textit{basic} cones
(constructed in Lemma \ref{LEMMA1}), we have
\begin{equation*}
\begin{array}{c}
\bigcup\limits_{\varepsilon_{i}\in\left\{ 0,1\right\} ,\ 1\leq i\leq d-1}%
\text{Gen}\left(
C_{\varepsilon_{1},\ldots,\varepsilon_{d-1}}^{\left(
d\right) }\right) =\mathcal{L}_{\mathbf{m}}^{\left( d\right) }\,.%
\end{array}
\end{equation*}
Hence, every element of $\tau_{P_{\mathbf{m}}^{\left( d\right) }}^{\vee}\cap(%
\mathbb{Z}^{d})^{\vee}$ can be written as non-negative integral
linear combination of the elements of
$\mathcal{L}_{\mathbf{m}}^{\left( d\right) }.$
\hfill$\square\bigskip$

\noindent Now let $\mathcal{A}_{\mathbf{m}}^{\left( d\right) }=\left( \text{I%
}_{d}\ \ \ \mathbf{M}^{\intercal}\right) $ be the $d\times\left(
2d-1\right) $ integral matrix with
\begin{equation*}
\mathbf{M:=}\left(
\begin{array}{ccccccc}
m_{1,1} & -1 & 0 & \cdots & \cdots & 0 & 0 \\
m_{2,1} & m_{2,2} & -1 & \cdots & \cdots & 0 & 0 \\
m_{3,1} & m_{3,2} & \ddots & \ddots &  & \vdots & \vdots \\
\vdots & \vdots & \ddots & \ddots & \ddots & \vdots & \vdots \\
\vdots & \vdots & \vdots & \ddots & \ddots & -1 & 0 \\
m_{d-1,1} & m_{d-1,2} & m_{d-1,3} & \cdots & m_{d-1,d-2} & m_{d-1,d-1} & -1%
\end{array}
\!\!\right) ,
\end{equation*}
and I$_{d}$ the $(d\times d)$-identity matrix. {}(Notice that the
transposes of the column vectors of
$\mathcal{A}_{\mathbf{m}}^{\left( d\right) \,}$ are precisely the
elements of the set $\mathcal{L}_{\mathbf{m}}^{\left( d\right)
}$). Consider the integer lattice
$\Lambda_{\mathcal{L}_{\mathbf{m}}^{\left( d\right) }}$ w.r.t.
$\mathcal{L}_{\mathbf{m}}^{\left( d\right) }$, i.e.,
\begin{align*}
\Lambda_{\mathcal{L}_{\mathbf{m}}^{\left( d\right) }} & =\left\{ \mathbf{\ell%
}=\left( \ell_{1},\ldots,\ell_{2d-1}\right) \in\mathbb{Z}^{2d-1}\
\left\vert \ \sum_{\kappa=1}^{d}\,\ell_{\kappa}\,e_{\kappa}^{\vee
}+\sum_{i=1}^{d-1}\,\ell_{d+i}\,\left( m_{i}-e_{i+1}^{\vee}\right)
\right.
=0\right\} \smallskip \\
& =\left\{ \mathbf{\ell}=\left( \ell_{1},\ldots,\ell_{2d-1}\right) \in%
\mathbb{Z}^{2d-1}\ \left\vert \ \mathcal{A}_{\mathbf{m}}^{\left( d\right) }\,%
\mathbf{\ell}^{\intercal}=0\right. \right\} =\text{Ker}\left(
\psi\right) ,
\end{align*}
where $\psi$ is the homomorphism: $\mathbb{Z}^{2d-1}\ni\mathbf{\ell\,}\longmapsto\psi\left( \mathbf{\ell }%
\right) =\,\mathcal{A}_{\mathbf{m}}^{\left( d\right) }\,\mathbf{\ell }%
^{\intercal}\in\mathbb{Z}^{d}.$
%\begin{equation*}
%\end{equation*}
By Lemma \ref{LEMMA2}, the characters $\mathbf{e}\left(
e_{1}^{\vee}\right),\ldots,\mathbf{e}\left( e_{d}^{\vee}\right)
,\,\mathbf{e}\left( m_{1}-e_{2}^{\vee}\right) ,\mathbf{e}\left(
m_{2}-e_{3}^{\vee}\right) ,\ldots,\mathbf{e}\left(
m_{d-1}-e_{d}^{\vee}\right)$ generate $\mathbb{C}[\tau_{P_{\mathbf{m}}^{\left( d\right) }}^{\vee}\cap(%
\mathbb{Z}^{d})^{\vee}]$. Hence, the affine toric variety
$U_{\tau_{P}}\cong
U_{\tau_{P_{\mathbf{m}}^{\left( d\right) }}}$ admits an embedding into $%
\mathbb{C}^{2d-1}$ w.r.t. $\mathcal{L}_{\mathbf{m}}^{\left(
d\right) },$ and
is, in particular, a g.c.i. of $d-1$ binomials (by Thm. \ref{Nak-thm}, Rem. %
\ref{REMNAK}(ii), and Thm. \ref{EMB}). The map
\begin{equation*}
\theta:\mathbb{C}\left[ z_{1},z_{2},\ldots,z_{d},z_{d+1},\ldots ,z_{2d-1}%
\right] \longrightarrow\mathbb{C}[\tau_{P_{\mathbf{m}}^{\left(
d\right) }}^{\vee}\cap(\mathbb{Z}^{d})^{\vee}]
\end{equation*}
defined by
\begin{equation*}
\theta\left( z_{\kappa}\right) :=\mathbf{e}\left(
e_{\kappa}^{\vee}\right) ,\ \ \forall\kappa,\ \
\kappa\in\{1,\ldots,d\},\text{\ }
\end{equation*}
and%
\begin{equation*}
\theta\left( z_{d+i}\right) :=\mathbf{e}\left(
m_{i}-e_{i+1}^{\vee}\right) ,\ \ \forall i,\ \
i\in\{1,\ldots,d-1\},
\end{equation*}
is a $\mathbb{C}$-algebra epimorphism. Let $\mathcal{I}:=\mathcal{I}_{%
\mathcal{A}_{\mathbf{m}}^{\left( d\right) }}:=$ Ker$\left(
\theta\right) $ denote its kernel. The column-vectors, say
$b_{1},b_{2},\ldots,b_{d-1}$, of
the $\left( 2d-1\right) \times\left( d-1\right) $-matrix%
\begin{equation*}
\mathcal{B}_{\mathbf{m}}^{\left( d\right) }:=\left(
\begin{array}{c}
\mathbf{M}^{\intercal} \\
-\text{ I}_{d-1}%
\end{array}
\right)=
\begin{pmatrix}
m_{1,1} & m_{2,1} & m_{3,1} & \dots & m_{d-2,1} & m_{d-1,1} \\
-1 & m_{2,2} & m_{3,2} & \dots & m_{d-2,2} & m_{d-1,2} \\
0 & -1 & m_{3,3} & \dots & m_{d-2,3} & m_{d-1,3} \\
0 & 0 & -1 & \ddots & m_{d-2,4} & m_{d-1,4} \\
\vdots & \vdots & \ddots & \ddots & \vdots & \vdots \\
0 & \ldots & 0 & 0 & -1 & m_{d-1,d-1} \\
0 & \ldots & \ldots & 0 & 0 & -1 \\
-1 & 0 & 0 & \dots & 0 & 0 \\
0 & -1 & 0 & \dots & 0 & 0 \\
0 & 0 & -1 & \ddots & 0 & 0 \\
\vdots & \vdots & \ddots & \ddots & \vdots & \vdots \\
0 & \ldots & 0 & 0 & -1 & 0 \\
0 & \ldots & \ldots & 0 & 0 & -1%
\end{pmatrix}
\end{equation*}
build up a $\mathbb{Z}$-basis of $\Lambda_{\mathcal{L}_{\mathbf{m%
}}^{\left( d\right) }}$. Let
$\mathcal{J}_{\mathcal{B}_{\mathbf{m}}^{\left( d\right) }}$ be the
lattice ideal of $\mathbb{C}\left[
z_{1},\ldots,z_{d},\ldots,z_{2d-1}\right] $ which is associated to $%
\mathcal{B}_{\mathbf{m}}^{\left( d\right) }$ (cf. \ref{LATIDEAL}).
In order to determine a generating system of $\mathcal{I}$
consisting of binomials whose exponents are expressed in terms of
the entries of our initial admissible sequence of free parameters
(\ref{NAKMATR}), it seems to be
reasonable to specify the saturation of $\mathcal{J}_{\mathcal{B}_{\mathbf{m}%
}^{\left( d\right) }}$ w.r.t. the product
$\prod_{j=1}^{2d-1}z_{j}\,$ of all available variables (see Thm.
\ref{SATUR}(i)). Nevertheless, this method would be rather
laborious from the computational point of view, because it
would involve elimination techniques or even primary decompositions of $%
\mathcal{J}_{\mathcal{B}_{\mathbf{m}}^{\left( d\right) }}$, and
would necessarily demand to perform a relatively high number of
Gr\"{o}bner basis algorithms (see \cite{BRS}, \cite{H-Sh},
\cite{Sturmfels1}). Instead, we
shall pass to another $\mathbb{Z}$-basis of $\Lambda_{\mathcal{L}_{\mathbf{m}%
}^{\left( d\right) }}$ \ whose matrix $\widehat{\mathcal{B}}_{\mathbf{m}%
}^{\left( d\right) }$ is dominating and we shall apply Thm.
\ref{SATUR}(ii).
(For affine semigroup rings which are complete intersections the \textit{%
existence} of an integral basis of their relation space having a
dominating coefficient matrix is guaranteed by a result of Fischer
and Shapiro, cf.
\cite[Cor. 2.10, p. 47]{F-Sh}. In the case at hand, $\widehat{\mathcal{B}}_{%
\mathbf{m}}^{\left( d\right) }$ will be constructed
\textit{explicitly}.)

\begin{lemma}
\label{LEMMALFABETAS}All integral $\left( 2d-1\right) \times\left(
d-1\right) $-matrices of the form \smallskip\
\begin{equation}
\left(
\begin{array}{cccccc}
\alpha_{1,1} & \alpha_{2,1} & \alpha_{3,1} & \cdots &
\alpha_{d-2,1} &
\alpha_{d-1,1} \\
-1 & \alpha_{2,2} & \alpha_{3,2} & \cdots & \alpha_{d-2,2} &
\alpha_{d-1,2}
\\
0 & -1 & \alpha_{3,3} & \cdots & \alpha_{d-2,3} & \alpha_{d-1,3} \\
0 & 0 & -1 & \ddots & \alpha_{d-2,4} & \alpha_{d-1,4} \\
\vdots & \vdots & \ddots & \ddots & \vdots & \vdots \\
0 & \cdots & 0 & 0 & -1 & \alpha_{d-1,d-1} \\
0 & \cdots & \cdots & 0 & 0 & -1 \\
-1 & \beta_{2,2} & \beta_{3,2} & \cdots & \beta_{d-2,2} & \beta_{d-1,2} \\
0 & -1 & \beta_{3,3} & \cdots & \beta_{d-2,3} & \beta_{d-1,3} \\
0 & 0 & -1 & \ddots & \beta_{d-2,4} & \beta_{d-1,4} \\
\vdots & \vdots & \ddots & \ddots & \vdots & \vdots \\
0 & \cdots & 0 & 0 & -1 & \beta_{d-1,d-1} \\
0 & \cdots & \cdots & 0 & 0 & -1%
\end{array}
\right) \,\smallskip   \label{ALFABETAS}
\end{equation}
where
\begin{equation*}
\alpha_{i,j}, \ \ 1\leq i\leq d-1, \ \ 1\leq j\leq i\leq d-1,
\end{equation*}
and
\begin{equation*}
\beta_{i,j}, \ \ 2\leq j\leq i\leq d-1,
\end{equation*}
are non-negative, and each of their columns contains at least one
positive element, are dominating matrices.
\end{lemma}

\noindent {}\textit{Proof}. Suppose that such an (obviously mixed)
matrix contains a mixed $\left( \rho \times \rho \right)
$-submatrix, for $\rho \in \{2,\ldots ,d-1\}$, with
\textquotedblleft column indices\textquotedblright\ $l_{1},\ldots
,l_{\rho }$, where
\begin{equation*}
l_{1}<l_{2}<\cdots <l_{\rho }.
\end{equation*}%
Then the negative entries of the column having index $l_{i}$ are
to be found in the rows whose indices belong to the set
\begin{equation*}
\mathfrak{N}_{i}:=\{l_{i}+1,l_{i}+d\},
\end{equation*}%
for every $i\in \{1,\ldots ,\rho \}.$ Since the $\left( \rho
\times \rho \right) $-submatrix under consideration is assumed to
be mixed, it has to contain in its first column a positive entry
which is located in the rows whose indices are within
\begin{equation*}
\mathfrak{N}_{\rho +1}:=\{1,\ldots ,l_{1},d+1,\ldots ,l_{1}+d-1\}.
\end{equation*}%
Since $l_{\rho }\leq d-1$, the $\rho +1$ sets $\left\{
\mathfrak{N}_{i}\ \left\vert \ 1\leq i\leq \rho +1\right. \right\}
$ are pairwise disjoint, but our $\left( \rho \times \rho \right)
$-submatrix must contain a \textquotedblleft row
index\textquotedblright\ from each of these sets,
which is impossible. Consequently, all integral matrices of the form (\ref%
{ALFABETAS}) are dominating matrices.\hfill $\square $

\begin{remark}
Our intention is to prove that after having performed (at most)
$d-2$
suitable unimodular transformations to the entries of our initial matrix $%
\mathcal{B}_{\mathbf{m}}^{\left( d\right) }$, we construct a $\mathbb{Z}$%
-basis of $\Lambda_{\mathcal{L}_{\mathbf{m}}^{\left( d\right) }}$
\ whose
matrix $\widehat{\mathcal{B}}_{\mathbf{m}}^{\left( d\right) }$ is of type (%
\ref{ALFABETAS}). This procedure will be realized in three steps.
In the first step, which explains where our motivation comes from,
we discuss what happens in the ``low'' dimensions $d=3$ and $d=4$.
In the second step, we present the recursive principle by means of
which we modify the last column of
$\mathcal{B}_{\mathbf{m}}^{\left( d\right) }$. Finally, in the
third step we apply unimodular transformations of the same sort to
the next coming column vectors. (In the particular case in which
\textit{all} the entries of (\ref{NAKMATR}) are non-negative,
$\mathcal{B}_{\mathbf{m}}^{\left( d\right)
}$ is itself dominating and there is no need to proceed, cf. \ref{REMMAIN}%
(iii)).
\end{remark}

\noindent$\blacktriangleright$ {}\textbf{Step 1: Low Dimensions.
}Let us start with $d=3$. In this case,
\begin{equation*}
\mathcal{B}_{\mathbf{m}}^{\left( 3\right) }=%
\begin{pmatrix}
m_{1,1} & m_{2,1} \\
-1 & m_{2,2} \\
0 & -1 \\
-1 & 0 \\
0 & -1%
\end{pmatrix}%
\begin{array}{l}
\\
\\
\\
\\
.%
\end{array}
\end{equation*}
If $m_{2,2}\geq0$ nothing is to do. So we may assume that
$m_{2,2}<0$. Adding $m_{2,2}$ times the first column to the second
one (which corresponds to a unimodular transformation) results in
\begin{equation*}
\begin{pmatrix}
m_{1,1} & m_{2,1}+m_{1,1}m_{2,2} \\
-1 & 0 \\
0 & -1 \\
-1 & -m_{2,2} \\
0 & -1%
\end{pmatrix}%
\begin{array}{l}
\\
\\
\\
\\
.%
\end{array}
\end{equation*}
In view of the ``admissibility conditions'' (\ref{ADM3}) the
matrix is of type (\ref{ALFABETAS}). (We should mention at this
point that, for $d=3$, Ishida makes similar choices by using some
purely geometric arguments, cf. \cite[proof of Thm. 8.1, in
particular pp. 140-141]{Ishida}).\smallskip

\noindent{}Let us now increase the dimension by one. For $d=4$,
\begin{equation*}
\mathcal{B}_{\mathbf{m}}^{\left( 4\right) }=%
\begin{pmatrix}
m_{1,1} & m_{2,1} & m_{3,1} \\
-1 & m_{2,2} & m_{3,2} \\
0 & -1 & m_{3,3} \\
0 & 0 & -1 \\
-1 & 0 & 0 \\
0 & -1 & 0 \\
0 & 0 & -1%
\end{pmatrix}%
\begin{array}{l}
\\
\\
\\
\\
\\
\\
.%
\end{array}
\end{equation*}
We start by looking at the last column. Assume, for instance, that $%
m_{3,3}<0 $. Adding $m_{3,3}$ times the second column to the last
one we obtain
\begin{equation*}
\begin{pmatrix}
m_{1,1} & m_{2,1} & m_{3,1}+m_{2,1}m_{3,3} \\
-1 & m_{2,2} & m_{3,2}+m_{2,2}m_{3,3} \\
0 & -1 & 0 \\
0 & 0 & -1 \\
-1 & 0 & 0 \\
0 & -1 & -m_{3,3} \\
0 & 0 & -1%
\end{pmatrix}%
\begin{array}{l}
\\
\\
\\
\\
\\
\\
.%
\end{array}
\end{equation*}
If
\begin{equation*}
m_{3,2}+m_{2,2}m_{3,3}\geq0,
\end{equation*}
the ``admissibility conditions'' (\ref%
{ADM4}) tell us that the entries of the last column of our matrix
are like
those of the last column of the matrices of type (\ref{ALFABETAS}). If
\begin{equation*}
m_{3,2}+m_{2,2}m_{3,3}<0,
\end{equation*}
then we add $\left( m_{3,2}+m_{2,2}m_{3,3}\right)$-times the first
row to the last one and we get
\begin{equation*}
\begin{pmatrix}
m_{1,1} & m_{2,1} &
m_{3,1}+m_{2,1}m_{3,3}+m_{1,1}(m_{3,2}+m_{2,2}m_{3,3})
\\
-1 & m_{2,2} & 0 \\
0 & -1 & 0 \\
0 & 0 & -1 \\
-1 & -0 & -(m_{3,2}+m_{2,2}m_{3,3}) \\
0 & -1 & -m_{3,3} \\
0 & 0 & -1%
\end{pmatrix}%
\begin{array}{l}
\\
\\
\\
\\
\\
\\
.%
\end{array}
\end{equation*}
Again by (\ref{ADM4}) the last column has the desired property.
After that we can start to transform the first two columns as we
did before in dimension $3$.

Next, let us assume that
\begin{equation*}
m_{3,3}\geq0, \ \ \text{but} \ \ m_{3,2}<0.
\end{equation*}
In this case, we add $m_{3,2}$ times the first column to the last
one and we get
\begin{equation*}
\begin{pmatrix}
m_{1,1} & m_{2,1} & m_{3,1}+m_{1,1}m_{3,2} \\
-1 & m_{2,2} & 0 \\
0 & -1 & m_{3,3} \\
0 & 0 & -1 \\
-1 & -0 & -m_{3,2} \\
0 & -1 & 0 \\
0 & 0 & -1%
\end{pmatrix}%
\begin{array}{l}
\\
\\
\\
\\
\\
\\
.%
\end{array}
\end{equation*}
Again the entries of the last column of our matrix are like those
of the last column of the matrices of type
(\ref{ALFABETAS}).\noindent\bigskip

\noindent{}{}{}$\blacktriangleright$ {}\textbf{Step 2: The
recursive principle. }We define appropriate matrix operations, so
that the last column looks like in (\ref{ALFABETAS}). Since these
operations (as we shall see below) do not affect the other
columns, we can apply a recursive argument.

More precisely, for $i=d-2,\dots,1$ we define recursively some
non-positive integers $\lambda_{i}$ and some non-negative integers
$\mu_{i}$ as follows:
\begin{equation*}
\begin{array}{l}
\lambda_{i}:=\left\{
\begin{array}{ll}
\min\left\{ 0,m_{d-1,d-1}\right\} , & \text{if }i=d-2, \\
\  & \  \\
\min\left\{
0,m_{d-1,i+1}+\sum\limits_{k=i+1}^{d-2}\lambda_{k}\,m_{k,i+1}\right\}
, &
\text{if }i\leq d-3,%
\end{array}
\right. \\
\  \\
\mu_{i}:=\left\{
\begin{array}{ll}
\max\left\{ 0,m_{d-1,d-1}\right\} , & \text{if }i=d-2, \\
\  & \  \\
\max\left\{
0,m_{d-1,i+1}+\sum\limits_{k=i+1}^{d-2}\lambda_{k}\,m_{k,i+1}\right\}
, &
\text{if }i\leq d-3.%
\end{array}
\right.%
\end{array}
\end{equation*}
Furthermore, we set
\begin{equation*}
\mu_{0}:=m_{d-1,1}+\sum_{i=1}^{d-2}\lambda_{i}\,m_{i,1}.
\end{equation*}
Since
\begin{equation*}
\mu_{i}=-\lambda_{i}+\left( m_{d-1,i+1}+\sum_{k=i+1}^{d-2}\lambda
_{k}\,m_{k,i+1}\right) ,\quad i=1,\dots,d-2,
\end{equation*}
the vector
$\widehat{b}_{d-1}:=b_{d-1}+\sum_{i=1}^{d-2}\lambda_{i}\,b_{i}$
can be written as
\begin{equation*}
\widehat{b}_{d-1}=\left(
\mu_{0},\mu_{1},\dots,\mu_{d-2},-1,-\lambda
_{1},\dots,-\lambda_{d-2},-1\right) ^{\intercal}.
\end{equation*}

\begin{lemma}
$\widehat{b}_{d-1}$ contains both negative and positive entries.
\end{lemma}

\noindent{}\textit{Proof}. Apparently, $\widehat{b}_{d-1}$ has two
negative coordinates. Suppose that $\widehat{b}_{d-1}$ has no
positive entries. Then
we have $\lambda_{i}=\mu_{i}=0$, $i=1,\dots,d-2$, which implies $%
m_{d-1,i+1}=0$, for $i=d-2,\dots,1$, and $\mu_{0}=m_{d-1,1}$.
However, by
the definition of the free parameters $\mathbf{m}$ (see Definition \ref%
{FREPAR}) we know that $m_{d-1}\neq0$ and from (\ref{ADMISSIBLE}) we get $%
m_{d-1,1}=\left\langle m_{d-1},e_{1}\right\rangle \geq0$, which
leads to a contradiction.\hfill $\square\medskip$

Next, we shall show that the first entry of $\widehat{b}_{d-1}$, i.e., $%
\mu_{0}$, is also non-negative. To this end we define
\begin{equation*}
\varepsilon_{i}:=\left\{
\begin{array}{lll}
0, & \text{ if }\lambda_{i}=0 &  \\
1, & \text{ if }\lambda_{i}<0 &
\end{array}
\right. ,\quad\text{\ for all }i\in\{1,\dots,d-2\},
\end{equation*}
and write $\lambda_{i}$ as follows:
\begin{equation*}
\lambda_{i}=\varepsilon_{i}\left(
m_{d-1,i+1}+\sum_{k=i+1}^{d-2}\lambda _{k}\,m_{k,i+1}\right) .
\end{equation*}
Now for $n\in\{1,\ldots,d-2\}$ and $k\in\{0,1,\ldots,d-2\}$ we set
\begin{equation*}
p_{k,n}:=\sum_{n=i_{0}<i_{1}<\cdots<i_{k}<i_{k+1}=d-1}\,\,\prod_{j=1}^{k+1}%
\varepsilon_{i_{j-1}}\,m_{i_{j},i_{j-1}+1}.
\end{equation*}
Note that for all $n\in\{1,\ldots,d-2\}$, we obtain%
\begin{equation*}
p_{0,n}=\varepsilon_{n}\cdot m_{d-1,n+1}.
\end{equation*}

\begin{lemma}
\label{LAMDASUM}For all $n\in\{1,2,\ldots,d-2\}$ we have%
\begin{equation*}
\lambda_{n}=\sum_{k=0}^{d-2-n}p_{k,n}.
\end{equation*}
\end{lemma}

\noindent{}\textit{Proof}. First we check the identity
\begin{equation}
p_{l+1,n-1}=\sum_{k=n}^{d-2-l}\varepsilon_{n-1}\,m_{k,n}\,p_{l,k}.
\label{IDENTP}
\end{equation}
Again this follows immediately from the definition, since
\begin{equation*}
\begin{split}
& \sum_{k=n}^{d-2-l}\varepsilon_{n-1}\,m_{k,n}\,p_{l,k}\smallskip \\
&
=\sum_{k=n}^{d-2-l}\varepsilon_{n-1}\,m_{k,n}\sum_{k=i_{0}<i_{1}<%
\cdots<i_{l}<i_{l+1}=d-1}\prod_{j=1}^{l+1}\varepsilon_{i_{j-1}}%
\,m_{i_{j},i_{j-1}+1}\smallskip \\
&
=\sum_{k=n}^{d-2-l}\varepsilon_{n-1}\,m_{k,n}\sum_{k=i_{1}<i_{2}<%
\cdots<i_{l+1}<i_{l+2}=d-1}\prod_{j=2}^{l+2}\varepsilon_{i_{j-1}}%
\,m_{i_{j},i_{j-1}+1}\smallskip \\
&
=\sum_{n-1=i_{0}<i_{1}<i_{2}<\cdots<i_{l+1}<i_{l+2}=d-1}\prod_{j=1}^{l+2}%
\varepsilon_{i_{j-1}}\,m_{i_{j},i_{j-1}+1}=p_{l+1,n-1}.
\end{split}
\end{equation*}
To prove the proposition we apply (backwards) induction with
respect to $n$. If $n=d-2,$ then we have
$\lambda_{d-2}=\varepsilon_{d-2}\,m_{d-1,d-1}$. So let $n$ be
$<d-2$. Then we may write
\begin{equation*}
\begin{split}
\lambda_{n-1} & =\varepsilon_{n-1}\left(
m_{d-1,n}+\sum_{k=n}^{d-2}\lambda_{k}\,m_{k,n}\right)=p_{0,n-1}+\sum_{k=n}^{d-2}\sum_{l=0}^{d-2-k}p_{l,k}\varepsilon
_{n-1}\,m_{k,n}\smallskip \\
&
=p_{0,n-1}+\sum_{l=0}^{d-2-n}\sum_{k=n}^{d-2-l}p_{l,k}\varepsilon
_{n-1}\,m_{k,n} \\
& =p_{0,n-1}+\sum_{l=0}^{d-2-n}p_{l+1,n-1},
\end{split}
\end{equation*}
where the last equation follows from
(\ref{IDENTP}).\hfill$\square$

\begin{lemma}
\begin{equation*}
\mu_{0}\geq0.
\end{equation*}
\label{mainlemma}
\end{lemma}

\noindent{}\textit{Proof}. From the definition of $\mu_{0}$ and Lemma \ref%
{LAMDASUM} we get
\begin{equation}
\begin{split}
\mu_{0} &
=m_{d-1,1}+\sum_{n=1}^{d-2}\lambda_{n}\,m_{n,1}=m_{d-1,1}+\sum_{n=1}^{d-2}%
\sum_{k=0}^{d-2-n}p_{k,n}\,m_{n,1}\smallskip \\
&
=m_{d-1,1}+\sum_{k=0}^{d-3}\sum_{n=1}^{d-2-k}p_{k,n}\,m_{n,1}=m_{d-1,1}+%
\sum_{k=1}^{d-2}\sum_{n=1}^{d-1-k}p_{k-1,n}\,m_{n,1}.
\end{split}
\label{FIRSTPART}
\end{equation}
Furthermore, by the definition of $p_{k,n}$ we have
\begin{equation*}
\begin{split}
\sum_{n=1}^{d-1-k}p_{k-1,n}\,m_{n,1} & =\sum_{n=1}^{d-1-k}m_{n,1}\hspace{%
-0.5cm}\sum_{n=i_{0}<i_{1}<\cdots<i_{k-1}<i_{k}=d-1}\,\,\prod
_{j=1}^{k}\varepsilon_{i_{j-1}}\,m_{i_{j},i_{j-1}+1}\smallskip \\
&
=\sum_{n=1}^{d-1-k}m_{n,1}\hspace{-0.5cm}\sum_{n=i_{1}<i_{2}<\cdots
<i_{k}<i_{k+1}=d-1}\,\,\prod_{j=2}^{k+1}\varepsilon_{i_{j-1}}%
\,m_{i_{j},i_{j-1}+1}\smallskip \\
&
=\sum_{0=i_{0}<i_{1}<i_{2}<\cdots<i_{k}<i_{k+1}=d-1}\,m_{i_{1},i_{0}+1}\,%
\prod_{j=2}^{k+1}\varepsilon_{i_{j-1}}\,m_{i_{j},i_{j-1}+1}\smallskip \\
&
=\sum_{0=i_{0}<i_{1}<i_{2}<\cdots<i_{k}<i_{k+1}=d-1}\,m_{d-1,i_{k}+1}\,%
\prod_{j=1}^{k}\varepsilon_{i_{j}}\,m_{i_{j},i_{j-1}+1}.
\end{split}
\end{equation*}
This, combined with (\ref{FIRSTPART}), gives
\begin{equation*}
\mu_{0}=\sum_{k=0}^{d-2}\quad\sum_{0=i_{0}<i_{1}<i_{2}<%
\cdots<i_{k}<i_{k+1}=d-1}m_{d-1,i_{k}+1}\cdot\quad\prod_{j=1}^{k}%
\varepsilon_{i_{j}}\,m_{i_{j},i_{j-1}+1},
\end{equation*}
which is non-negative by Corollary
\ref{POSITIVITY}.\hfill$\square\bigskip$

What we have done so far can be summarized in the following:

\begin{corollary}
\label{BUCOR}Let $\mathcal{U}^{\left( d\right) }\in$ \emph{GL}$\left( d-1,%
\mathbb{Z}\right) $ be the unimodular matrix given by%
\begin{equation*}
\mathcal{U}^{\left( d\right) }=\left(
\begin{array}{ccccc}
1 & 0 & 0 & \cdots & \lambda_{1} \\
0 & 1 & 0 & \cdots & \lambda_{2} \\
\vdots & \ddots & \ddots & \vdots & \vdots \\
0 & \cdots & 0 & 1 & \lambda_{d-2} \\
0 & 0 & \cdots & 0 & 1%
\end{array}
\right)
\begin{array}{l}
\\
\\
\\
\\
.%
\end{array}
\end{equation*}
Then we have%
\begin{equation*}
\mathcal{B}_{\mathbf{m}}^{\left( d\right) }\mathcal{U}^{\left(
d\right) }=\left(
\begin{array}{cccccc}
m_{1,1} & m_{2,1} & m_{3,1} & \cdots & m_{d-2,1} & \mu_{0} \\
-1 & m_{2,2} & m_{3,2} & \cdots & m_{d-2,2} & \mu_{1} \\
0 & -1 & m_{3,3} & \cdots & m_{d-2,3} & \mu_{2} \\
0 & 0 & -1 & \ddots & m_{d-2,4} & \mu_{3} \\
\vdots & \vdots & \ddots & \ddots & \vdots & \vdots \\
0 & \cdots & 0 & 0 & -1 & \mu_{d-2} \\
0 & \cdots & \cdots & 0 & 0 & -1 \\
-1 & 0 & 0 & \cdots & 0 & -\lambda_{1} \\
0 & -1 & 0 & \cdots & 0 & -\lambda_{2} \\
0 & 0 & -1 & \ddots & 0 & -\lambda_{3} \\
\vdots & \vdots & \ddots & \ddots & \vdots & \vdots \\
0 & \cdots & 0 & 0 & -1 & -\lambda_{d-2} \\
0 & \cdots & \cdots & 0 & 0 & -1%
\end{array}
\right) \ \smallskip
\end{equation*}
where $\mu_{i}\geq0$ for all $i\in\{0,1,\ldots,d-2\}$,
$\lambda_{i}\leq0$ for all $i\in\{1,\ldots,d-2\}$, and the last
column contains both positive and negative entries.
\end{corollary}

We observe that the first $d-2$ columns of the matrix $\mathcal{B}_{\mathbf{m%
}}^{\left( d\right) }\mathcal{U}^{\left( d\right) }$ are the same
as the
columns of the matrix $\mathcal{B}_{\mathbf{m}}^{\left( d\right) }$, namely%
\begin{equation*}
\begin{array}{ccccc}
\left(
\begin{array}{c}
m_{1,1} \\
-1 \\
0 \\
0 \\
\vdots  \\
0 \\
0 \\
-1 \\
0 \\
0 \\
\vdots  \\
0 \\
0%
\end{array}%
\right) , & \left(
\begin{array}{c}
m_{2,1} \\
m_{2,2} \\
-1 \\
0 \\
\vdots  \\
0 \\
0 \\
0 \\
-1 \\
0 \\
\vdots  \\
0 \\
0%
\end{array}%
\right) , & \left(
\begin{array}{c}
m_{3,1} \\
m_{3,2} \\
m_{3,3} \\
-1 \\
\vdots  \\
0 \\
0 \\
0 \\
0 \\
-1 \\
\vdots  \\
0 \\
0%
\end{array}%
\right) , & \ldots \ , & \left(
\begin{array}{c}
m_{d-2,1} \\
m_{d-2,2} \\
m_{d-2,3} \\
m_{d-2,4} \\
\vdots  \\
-1 \\
0 \\
0 \\
0 \\
0 \\
\vdots  \\
-1 \\
0%
\end{array}%
\right)
\end{array}%
,
\end{equation*}%
and that the $d$-th and $\left( 2d-1\right) $-th row contain
non-zero entries only in the last column of
$\mathcal{B}_{\mathbf{m}}^{\left( d\right) }\mathcal{U}^{\left(
d\right) }$. This allows us to apply our transformations, which
have been carried out so far only for the last
column, successively to the other columns too.\bigskip

\noindent{}$\blacktriangleright$ {}{}\textbf{Step 3: Generalizing
the recursion for all column vectors. }It is enough to equip our
lambdas and mus with an additional index, just for keeping track
of the next coming columns
(viewed backwards). That's why we define recursively non-negative integers $%
\lambda_{i,j}$ and non-positive integers $\mu_{i,j}$ by the formulae (\ref%
{LAMDAS}) and (\ref{MIOUS}), respectively.
($\lambda_{i,d-1},\mu_{i,d-1}$ are exactly the numbers which we
called before $\lambda_{i}$ and $\mu_{i}$).
Moreover, using the column vectors $b_{1},\ldots,b_{d-1}$ of $\mathcal{B}_{%
\mathbf{m}}^{\left( d\right) }$, we introduce the integer linear
combinations:%
\begin{equation*}
\widehat{b}_{j}:=\left\{
\begin{array}{ll}
b_{1}, & \text{if }j=1 \\
b_{j}+\sum\limits_{i=1}^{j-1}\lambda_{i,j}\,b_{i}, & \text{if \ }%
j\in\{2,\ldots,d-1\}.%
\end{array}
\right.
\end{equation*}
Finally, we define $\mathcal{U}_{j}^{\left( d\right) }\in $
GL$\left( d-1,\mathbb{Z}\right) $ for $j=2,\dots,d-1$ as follows:
\begin{equation*}
\mathcal{U}_{j}^{\left( d\right) }:=%
\begin{pmatrix}
1 & 0 & 0 & \cdots & 0 & \lambda_{1,j} & 0 & \cdots & 0 \\
0 & 1 & 0 & \cdots & 0 & \lambda_{2,j} & 0 & \cdots & 0 \\
\vdots &  & \ddots &  & \vdots & \vdots & 0 & \cdots & 0 \\
\vdots &  & 0 & \ddots & \vdots & \vdots & 0 & \cdots & 0 \\
&  &  &  & \ddots & \lambda_{j-1,j} & 0 & \cdots & 0 \\
0 & 0 & 0 & \cdots & 0 & 1 & 0 & \cdots & 0 \\
\vdots &  &  &  &  &  & \ddots &  &  \\
&  &  &  &  &  &  & \ddots &  \\
0 & 0 & 0 & \cdots & \cdots & \cdots & \cdots & 0 & 1%
\end{pmatrix}
\end{equation*}
Each $\mathcal{U}_{j}^{\left( d\right) }$ is obviously
an upper triangular matrix, with $1$'s as diagonal elements$.$ All
the other non-trivial elements are contained in the $j$-th column.

\begin{proposition}
\label{DOMBDACH}Using the above matrices, the product%
\begin{equation*}
\widehat{\mathcal{B}}_{\mathbf{m}}^{\left( d\right) }:=\left( \widehat {b}%
_{1},\ldots,\widehat{b}_{d-1}\right)
=\mathcal{B}_{\mathbf{m}}^{\left(
d\right) }\cdot\mathcal{U}_{d-1}^{\left( d\right) }\cdot\mathcal{U}%
_{d-2}^{\left( d\right)
}\cdot\,\cdots\,\cdot\mathcal{U}_{2}^{\left( d\right) }
\end{equation*}
reads as%
\begin{equation}
\widehat{\mathcal{B}}_{\mathbf{m}}^{\left( d\right) }=\left(
\begin{array}{cccccc}
m_{1,1} & \mu_{0,2} & \mu_{0,3} & \cdots & \mu_{0,d-2} & \mu_{0,d-1} \\
-1 & \mu_{1,2} & \mu_{1,3} & \cdots & \mu_{1,d-2} & \mu_{1,d-1} \\
0 & -1 & \mu_{2,3} & \cdots & \mu_{2,d-2} & \mu_{2,d-1} \\
0 & 0 & -1 & \ddots & \mu_{3,d-2} & \mu_{3,d-1} \\
\vdots & \vdots & \ddots & \ddots & \vdots & \vdots \\
0 & \cdots & 0 & 0 & -1 & \mu_{d-2,d-1} \\
0 & \cdots & \cdots & 0 & 0 & -1 \\
-1 & -\lambda_{1,2} & -\lambda_{1,3} & \cdots & -\lambda_{1,d-2} &
-\lambda_{1,d-1} \\
0 & -1 & -\lambda_{3,3} & \cdots & -\lambda_{3,d-2} & -\lambda_{2,d-1} \\
0 & 0 & -1 & \ddots & -\lambda_{4,d-2} & -\lambda_{3,d-1} \\
\vdots & \vdots & \ddots & \ddots & \vdots & \vdots \\
0 & \cdots & 0 & 0 & -1 & -\lambda_{d-2,d-1} \\
0 & \cdots & \cdots & 0 & 0 & -1%
\end{array}
\right) \,\smallskip   \label{BDACHMATRIX}
\end{equation}
which is a dominating matrix with only $-1$'s as negative entries.
\end{proposition}

\noindent{}\textit{Proof}. According to the ``admissibility conditions'' (%
\ref{ADMISSIBLE}) we have always $m_{1,1}>0$. We observe that the
unimodular matrix $\mathcal{U}_{j}^{\left( d\right) }$ affects
only the $j$-th column; and since
$\mathcal{B}_{\mathbf{m}}^{\left( d\right) }$ is partitioned into
two upper triangular matrices, $\mathcal{U}_{j}^{\left( d\right)
}$affects
only the non-trivial elements of the $j$-th column of the matrix $\mathcal{B}%
_{\mathbf{m}}^{\left( d\right) }\cdot\mathcal{U}_{d-1}^{\left(
d\right) }\cdot\,\cdots\,\cdot\mathcal{U}_{j+1}^{\left( d\right)
}$. Hence, it suffices to use induction w.r.t. $d$ (by exploiting
Corollary \ref{BUCOR}),
and to take into account Lemma \ref{LEMMALFABETAS} and Lemma \ref{mainlemma}%
.\hfill$\square\bigskip$

\noindent {}\textbf{Proof of Theorem \ref{MAIN}}: By Proposition \ref%
{DOMBDACH}, $\widehat{\mathcal{B}}_{\mathbf{m}}^{\left( d\right)
}$ is a
dominating matrix and its column vectors $\widehat{b}_{1},\ldots ,\widehat{b}%
_{d-1}$ constitute a $\mathbb{Z}$-basis of $\Lambda _{\mathcal{L}_{\mathbf{m}%
}^{\left( d\right) }}$. Applying Theorem \ref{SATUR}(ii) we obtain
\begin{equation*}
\mathcal{I}=\mathcal{J}_{\widehat{\mathcal{B}}_{\mathbf{m}}^{\left(
d\right) }},
\end{equation*}%
and the $j$-th binomial of the constructed $d-1$ generators of
$\mathcal{I}$ is exactly that one containing the non-negative
entries of $\widehat{b}_{j}$ as exponents of the variables
$z_{1},\ldots ,z_{2d-1}$ in its first monomial and the opposites
of the non-positive entries of $\widehat{b}_{j}$ as exponents of
$z_{1},\ldots ,z_{2d-1}$ in its second monomial, $\forall
j,1\leq j\leq d-1,$cf. \ref{LATIDEAL} and (\ref{BDACHMATRIX}). Since $%
\widehat{\mathcal{B}}_{\mathbf{m}}^{\left( d\right) }$ has only
two $-1$'s as negative entries in each column, the second
monomials contain just two
variables, namely $z_{j+1}$ and $z_{d+j}$ for all $j\in \{1,\ldots ,d-1\}$%
.\hfill $\square $

\begin{remark}
By Lemma \ref{LEMMA2}, we have $\mathbf{Hilb}_{(\mathbb{Z}%
^{d})^{\vee}}(\tau_{P_{\mathbf{m}}^{\left( d\right) }}^{\vee})\subseteq%
\mathcal{L}_{\mathbf{m}}^{\left( d\right) },$ and the inclusion
may be strict. Next Lemma gives an explicit characterizarion of
the Hilbert basis
elements and allows us to realize the ``minimal'' embedding of $U_{\tau_{P_{%
\mathbf{m}}^{\left( d\right) }}}$ by eliminating redundant
variables.
\end{remark}

\begin{lemma}
\label{HILBIS}Let $v\in\mathcal{L}_{\mathbf{m}}^{\left( d\right)
}.$ Then
\begin{equation*}
v\in\mathbf{Hilb}_{(\mathbb{Z}^{d})^{\vee}}(\tau_{P_{\mathbf{m}}^{\left(
d\right) }}^{\vee})\Longleftrightarrow\left[ \nexists i\in\{2,\ldots ,d-1\}%
\text{\emph{, \ such that} }v=m_{i}\right]
\end{equation*}
\end{lemma}

\noindent {}\textit{Proof}. Suppose that $v=m_{i}$ for an index
$i\in \{2,\ldots ,d-1\}.$ Then $v$ can be written as the sum
\begin{equation*}
v=\left( m_{i}-e_{i+1}^{\vee }\right) +e_{i+1}^{\vee }
\end{equation*}%
which means that $v\notin \mathbf{Hilb}_{(\mathbb{Z}^{d})^{\vee }}(\tau _{P_{%
\mathbf{m}}^{\left( d\right) }}^{\vee })$ by (\ref{Hilbbasis}).
Conversely, assume that
\begin{equation*}
v\notin \mathbf{Hilb}_{(\mathbb{Z}^{d})^{\vee }}(\tau _{P_{\mathbf{m}%
}^{\left( d\right) }}^{\vee }).
\end{equation*}%
It is worth mentioning that $e_{d}^{\vee }$ and
$m_{d-1}-e_{d}^{\vee }$ are elements of
$\mathcal{L}_{\mathbf{m}}^{\left( d\right) }$ belonging always to
the Hilbert basis. We shall use again induction w.r.t. $d$. For
$d=3$ the assertion can be verified easily (cf. \cite[Prop. 8.7,
p. 138]{Ishida}). Now
let $d$ be $>3.$ Since $v\in \mathcal{L}_{\mathbf{m}}^{\left( d\right) }%
\mathbb{r}\mathbf{Hilb}_{(\mathbb{Z}^{d})^{\vee }}(\tau _{P_{\mathbf{m}%
}^{\left( d\right) }}^{\vee }),$ by Lemma \ref{LEMMA2} we may
express $v$ as a linear combination of the form
\begin{equation*}
v=\left( \sum_{\overline{v}\in \mathcal{L}_{\mathbf{m}}^{\left( d\right) }%
\mathbb{r}\{v,e_{d}^{\vee },m_{d-1}-e_{d}^{\vee }\}}\alpha _{\overline{v}}\,%
\overline{v}\right) +\beta \,e_{d}^{\vee }+\beta ^{\prime
}\,\left( m_{d-1}-e_{d}^{\vee }\right) ,
\end{equation*}%
for some $\alpha _{\overline{v}},\beta ,\beta ^{\prime }\in
\mathbb{Z}_{\geq 0}.$ Since $v\notin \{e_{d}^{\vee
},m_{d-1}-e_{d}^{\vee }\},$ its last coordinate has to be zero.
Therefore $\beta =\beta ^{\prime }$ with
\begin{equation}
v=\left( \sum_{\overline{v}\in \mathcal{L}_{\mathbf{m}}^{\left( d\right) }%
\mathbb{r}\{v,e_{d}^{\vee },m_{d-1}-e_{d}^{\vee }\}}\alpha _{\overline{v}}\,%
\overline{v}\right) +\beta \,m_{d-1}.  \label{LASTEQ}
\end{equation}%
The admissibility of $\mathbf{m}$ has as direct consequence that
$m_{d-1}\in \tau _{P_{\widetilde{\mathbf{m}}}^{\left( d-1\right)
}}^{\vee }$. Hence we have $v\in \tau
_{P_{\widetilde{\mathbf{m}}}^{\left( d-1\right) }}^{\vee }$, where
$\widetilde{\mathbf{m}}$ is the matrix $\mathbf{m}$ without the
last row.

On the other hand, by induction hypothesis, we may assume that
\begin{equation*}
v\in \mathbf{Hilb}_{(\mathbb{Z}^{d})^{\vee }}(\tau _{P_{\widetilde{\mathbf{m}%
}}^{\left( d-1\right) }}^{\vee })
\end{equation*}%
and hence, $v$ cannot be written as a non-trivial sum of elements
of $\tau
_{P_{\widetilde{\mathbf{m}}}^{\left( d-1\right) }}^{\vee }\cap (\mathbb{Z}%
^{d})^{\vee }$. So (\ref{LASTEQ}) is to be simplified as follows:
\begin{equation*}
v=\beta \,m_{d-1}.
\end{equation*}%
Since some coordinate of $v$ necessarily equals either $1$ or $-1,$ and $%
m_{d-1}\in (\mathbb{Z}^{d})^{\vee },$ we conclude that $\beta =1$,
as required.\hfill $\square $

\begin{examples}
Computing (by \ref{HILBIS}) the Hilbert bases of the duals of the
cones which support the Nakajima polytopes \ref{BEISPIELE} (i)
-(v) w.r.t. the rectangular lattice, we obtain:\smallskip\
\newline $\bullet$ For (i):
\begin{equation*}
\mathbf{Hilb}_{(\mathbb{Z}^{3})^{\vee}}(\tau_{P_{\mathbf{m}}^{\left(
3\right) }}^{\vee})=\mathcal{L}_{\mathbf{m}}^{\left( 3\right)
}\smallskip=\left\{
e_{1}^{\vee},e_{2}^{\vee},e_{3}^{\vee},2e_{1}^{\vee}-e_{2}^{\vee},2e_{1}^{%
\vee}+e_{2}^{\vee}-e_{3}^{\vee}\right\} .
\end{equation*}
$\bullet$ For (ii):\newline
\begin{equation*}
\mathbf{Hilb}_{(\mathbb{Z}^{3})^{\vee}}(\tau_{P_{\mathbf{m}}^{\left(
3\right) }}^{\vee})=\mathcal{L}_{\mathbf{m}}^{\left( 3\right) }\mathbb{r}%
\left\{ ke_{1}^{\vee}-e_{2}^{\vee\intercal}\right\} =\left\{
e_{1}^{\vee},e_{2}^{\vee},e_{3}^{\vee},ke_{1}^{\vee}-e_{2}^{\vee}-e_{3}^{%
\vee }\right\} .
\end{equation*}
$\bullet$ For (iii):%
\begin{equation*}
\mathbf{Hilb}_{(\mathbb{Z}^{4})^{\vee}}(\tau_{P_{\mathbf{m}}^{\left(
4\right) }}^{\vee})=\mathcal{L}_{\mathbf{m}}^{\left( 4\right) }\mathbb{r}%
\{e_{1}^{\vee}\},
\end{equation*}
\newline
where
\begin{equation*}
\mathcal{L}_{\mathbf{m}}^{\left( 4\right) }=\left\{
e_{1}^{\vee},e_{2}^{\vee},e_{3}^{\vee},e_{4}^{\vee},e_{1}^{\vee}-e_{2}^{%
\vee},e_{1}^{\vee
}-e_{3}^{\vee},2e_{1}^{\vee}-e_{2}^{\vee}-e_{3}^{\vee}-e_{4}^{\vee}\right\}
.
\end{equation*}
$\bullet$ For (iv):%
\begin{equation*}
\mathcal{L}_{\mathbf{m}}^{\left( d\right) }\smallskip=\{\left.
e_{i}^{\vee }\,\right\vert \ 1\leq i\leq
d\}\cup\{ke_{1}^{\vee}-e_{2}^{\vee}\}\cup\{\left.
e_{i}^{\vee}-e_{i+1}^{\vee}\ \right\vert \ 2\leq i\leq d-1\},
\end{equation*}
while%
\begin{equation*}
\mathbf{Hilb}_{(\mathbb{Z}^{d})^{\vee}}(\tau_{\mathbf{s}_{k}^{\left(
d\right)
}}^{\vee})=\{e_{1}^{\vee},e_{d}^{\vee}\}\cup\{ke_{1}^{\vee}-e_{2}^{\vee}\}%
\cup\{\left. e_{i}^{\vee}-e_{i+1}^{\vee}\ \right\vert \ 2\leq
i\leq d-1\}.
\end{equation*}
$\bullet$ For (v):%
\begin{equation*}
\begin{array}{l}
\mathbf{Hilb}_{(\mathbb{Z}^{d})^{\vee}}(\tau_{\mathbf{RP}\left(
k_{1},..,k_{d-1}\right) }^{\vee})=\mathcal{L}_{\mathbf{m}}^{\left(
d\right)
}\smallskip \\
=\{\left. e_{i}^{\vee}\ \right\vert \ 1\leq i\leq
d\}\cup\{k_{i}\left.
e_{i}^{\vee}-e_{i+1}^{\vee}\ \right\vert \ 1\leq i\leq d\}.\medskip%
\end{array}
\ \medskip
\end{equation*}
\end{examples}

\noindent {}\textbf{Proof of Corollary \ref{KOROLLAR}}: We define
$\mathfrak{Q}_{\mathbf{m}}^{\left( d\right) }$ and
$\mathfrak{R}_{\mathbf{m}}^{\left( d\right) }$ as in
(\ref{QINDEX}) and (\ref{RINDEX}), respectively. Suppose that one
of the unit vectors, say $e_{k}^{\vee },$ $k\in \left\{ 1,\ldots
,d-1\right\} ,$ does not belong to the Hilbert basis. Then by Lemma \ref%
{HILBIS} for some index $\gamma =\gamma _{k}$, which is $\geq k,$
we have $\ $necessarily $e_{k}^{\vee }=m_{\gamma }$, i.e., $k\in $
$\mathfrak{Q}_{\mathbf{m}}^{\left( d\right) }$. Looking at the
definition of $\lambda _{i,\gamma }$ and $\mu _{i,\gamma }$ we
find
\begin{equation*}
\lambda _{i,\gamma }=0\text{, \ \ \ for all \ }i\in \{0,\ldots
,\gamma -1\},
\end{equation*}%
and%
\begin{equation*}
\mu _{i,\gamma }=\left\{
\begin{array}{ll}
0, & \text{if }i\in \{0,\ldots ,\gamma -1\}\mathbb{r}\{k-1\}, \\
1, & \text{if }i=k-1,%
\end{array}%
\right.
\end{equation*}%
and the binomial corresponding to the $\gamma $-th column of $\widehat{%
\mathcal{B}}_{\mathbf{m}}^{\left( d\right) }$ equals
$z_{k}-z_{\gamma +1}z_{d+}{}_{\gamma }$; \ \ and conversely, if
one of the $d-1$ initial
binomials is of this type, then $e_{k}^{\vee }=m_{\gamma }$ with $%
e_{k}^{\vee }$ not belonging to the Hilbert basis. Analogously, if
a vector of type $m_{l}-e_{l+1}^{\vee }$, for some $l\in
\{1,\ldots ,d-2\}$, does not
belong to the Hilbert basis, then there exists an index $\delta =\delta _{l}$%
, which is $\geq l+1,$ so that $m_{l}-e_{l+1}^{\vee }=m_{\delta }$, i.e., $%
l\in $ $\mathfrak{R}_{\mathbf{m}}^{\left( d\right) }$. Again by the definition of $%
\lambda _{i,\delta }$ and $\mu _{i,\delta }$ we find%
\begin{equation*}
\lambda _{i,\delta }=\left\{
\begin{array}{ll}
0, & \text{if }i\in \{0,\ldots ,\delta -1\}\mathbb{r}\{l\}, \\
-1, & \text{if }i=l,%
\end{array}%
\right.
\end{equation*}%
and%
\begin{equation*}
\mu _{i,\delta }=0\text{, \ \ \ for all \ }i\in \{0,\ldots ,\delta
-1\},
\end{equation*}%
and the binomial corresponding to the $\delta $-th column of $\widehat{%
\mathcal{B}}_{\mathbf{m}}^{\left( d\right) }$ equals
$z_{d+l}-z_{\delta +1}z_{d+\delta }$; \ and conversely, if one of
the $d-1$ initial binomials
is of this type, then $m_{l}-e_{l+1}^{\vee }=m_{\delta }$ with $%
m_{l}-e_{l+1}^{\vee }$ not belonging to the Hilbert basis.
Obviously,
\begin{equation*}
\#(\mathbf{Hilb}_{(\mathbb{Z}^{d})^{\vee }}(\tau
_{P_{\mathbf{m}}^{\left( d\right) }}^{\vee
}))=\#(\mathcal{L}_{\mathbf{m}}^{\left( d\right)
})-\#\left( \mathfrak{Q}_{\mathbf{m}}^{\left( d\right) }\right) -\#(%
\mathfrak{R}_{\mathbf{m}}^{\left( d\right) }),
\end{equation*}%
and the assertion is true.\hfill $\square \bigskip $

\noindent{}\textbf{Comment and open problem}. The first partial
verification of the fact that the affine semigroup rings which are
complete intersections admit an ``inductive characterization''
\textit{in all dimensions} appeared already in the 1980's, in the
works of Watanabe \cite{Watanabe} and Nakajima \cite{Nakajima}
(who classified those which are invariant subrings of finite
abelian groups, and affine torus embeddings, respectively). This
was completely proved in 1997 by Fischer, Morris and Shapiro
\cite{F-M-Sh} via the theory of dominating matrices.

``Watanabe simplices'' (introduced in \cite{DHZ}) and, more
general, but in
a slightly different context, ``Nakajima polytopes'' provide a \textit{%
geometric parametrization} of the classes of semigroup rings treated in \cite%
{Watanabe} and \cite{Nakajima}, respectively. In the present
paper, based on Nakajima's classification, we ascertained that,
besides the above mentioned inductive characterization, there is
also some kind of \textit{recursion principle} governing a natural
set of generators of the relation space of the involved
semigroups. This gives rise to ask if this property can be
generalized for a wider class of affine semigroup rings, probably
in connection with a suitable combinatorial parametrization
resulting from partitions of minimal generating sets (also called
``semigroup gluings''), decomposition theorems of dominating
matrices or even from graph-theoretic
objects (like coloured paths etc), cf. \cite{R-G}, \cite{ROSALES}, \cite%
{F-M-Sh}, \cite{BMT}, and \cite{S-S-S}.%
\bigskip

%\noindent{}\textbf{Acknowledgements}. The second-named author would like to
%thank the Mathematics Department of the University of Crete for hospitality
%and support during the spring term 2001, where this work was initiated.

\end{document}